\documentclass[12pt]{amsart}

\usepackage{amsmath}
\usepackage{amssymb}

\allowdisplaybreaks                     
\theoremstyle{plain}
\newtheorem{theo}{Theorem}[section]
\newtheorem*{theorem*}{Theorem}
\newtheorem{lemm}[theo]{Lemma}

\newtheorem{prop}[theo]{Proposition}
\newtheorem{defi}[theo]{Definition}
\theoremstyle{remark}

\numberwithin{equation}{section}
\newcommand{\quadro}{\hfill\vrule height .9ex width .8ex depth -.1ex}

\newcommand{\NN}{\Bbb{N}}
\newcommand{\ZZ}{\Bbb{Z}}
\newcommand{\RR}{\Bbb{R}}
\newcommand{\CC}{\Bbb{C}}

\newcommand{\CZ}{Calder\'on--Zygmund }

\newcommand{\di}{\,{\rm{d}}}
\newcommand{\dir}{\,{\rm{d}}\rho}  
\newcommand{\dil}{\,{\rm{d}}\lambda}  

\newcommand{\nep}{{\rm{e}}} 

\begin{document}
\title[Boundedness from $H^1$ to $L^1$ of Riesz transforms]{\bf{Boundedness from $H^1$ to $L^1$ of Riesz transforms on a Lie group of exponential growth}}

\addtocounter{footnote}{1}
\addtocounter{footnote}{1}

\thanks{{\bf Acknowledgement.}  Work partially supported by the European Commission via the Network HARP, ``Harmonic analysis and related problems''. }

\author[P. Sj\"ogren and M. Vallarino]
{Peter Sj\"ogren and Maria Vallarino}

\maketitle

\begin{abstract}
Let $G$ be the Lie group $\RR^2\ltimes \RR^+$ endowed with the Riemannian symmetric space structure. Let $X_0,\, X_1,\,X_2$ be a distinguished basis of left-invariant vector fields of the Lie algebra of $G$ and define the Laplacian $\Delta=-(X_0^2+X_1^2+X_2^2)$. In this paper we consider the first order Riesz transforms $R_i=X_i\Delta^{-1/2}$ and $S_i=\Delta^{-1/2}X_i$, for $i=0,1,2$. We prove that the operators $R_i$, but not the $S_i$, are bounded from the Hardy space $H^1$ to $L^1$. We also show that the second-order Riesz transforms $T_{ij}=X_i\Delta^{-1}X_j$ are bounded from $H^1$ to $L^1$, while the Riesz transforms $S_{ij}=\Delta^{-1}X_iX_j$ and $R_{ij}=X_iX_j\Delta^{-1}$, for $i,j=0,1,2$, are not.
\end{abstract}

\maketitle

\section{Introduction}\label{intro}
Let $G$ be the Lie group $\RR^2\ltimes \RR^+$ where the product rule is the following:
$$(x_1,x_2,a)\cdot(x'_1,x'_2,a')=(x_1+a\,x'_1,x_2+a\,x'_2,a\,a') $$
for $(x_1,x_2,a),\,(x'_1,x'_2,a')\in G$. The group $G$ is not unimodular; the right and left Haar measures are given by
$$\dir(x_1,x_2,a)=a^{-1}\,\di x_1\di x_2\di a\qquad{\rm{and}}\qquad \dil(x_1,x_2,a)=a^{-3}\,\di x_1\di x_2\di a\,,$$
respectively. The modular function is thus $\delta(x_1,x_2,a)=a^{-2}$. Throughout this paper, unless explicitly stated, we consider the right measure $\rho$ on $G$ and we denote by $L^p$, $\|\cdot\|_p$ and $\langle\cdot,\cdot \rangle$  the $L^p$-space, the $L^p$-norm and the $L^2$-scalar product with respect to the measure $\rho$. 

The group $G$ has a Riemannian symmetric space structure, and the corresponding metric, which we denote by $d$, is that of the three-dimensional hyperbolic half-space. The metric $d$ is invariant under left translation and it is given by
\begin{equation}\label{metrica}
\cosh r(x_1,x_2,a)=\frac{a+a^{-1}+a^{-1}(x_1^2+x_2^2)}{2}\qquad \forall (x_1,x_2,a)\in G\,,
\end{equation}
where $r(x_1,x_2,a)=d\big((x_1,x_2,a),e  \big)$ denotes the distance of the point $(x_1,x_2,a)$ from the identity $e=(0,0,1)$ of $G$. It is easy to verify that if $r(x_1,x_2,a)<1$, then $r(x_1,x_2,a)\sim  |(x_1,x_2,\log a) |$, where $|\cdot |$ denotes the euclidean norm in $\RR^3$. The measure of a hyperbolic ball $B_r$, centered at the identity and of radius $r$, behaves like
$$\lambda(B_r)=\rho(B_r)\sim  \begin{cases}
r^3&{\rm{if}}~r<1\\
\nep^{2r}&{\rm{if}}~r\geq 1\,.
\end{cases}
$$
Thus $G$ is a group of exponential growth. In this context, the classical \CZ theory and the classical definition of the atomic Hardy space $H^1$ (see \cite{CW, St}) do not apply. Recently W.~Hebisch and T.~Steger \cite{HS} constructed a new \CZ theory which holds in some spaces of exponential growth, in particular in the space $(G,d,\rho)$ defined above. The main idea is to replace the family of balls which is used in the classical \CZ theory by a suitable family of parallelepipeds which we call \emph{\CZ sets}. The definition appears in \cite{HS} and implicitly in \cite{GS}, and reads as follows.
\begin{defi}\label{Czsets}
A Calder\'on-Zygmund set is a parallelepiped $R=[b_1-L/2,b_1+L/2]\times[b_2-L/2,b_2+L/2]\times [a\nep^{-r},a\nep^r]$, where $L>0$, $r>0$ and $a\in\RR^+$ are related by
$$\nep^{2}a\,r\leq L< \nep^{8 }a\,r\qquad{\rm{if~}}r<1\,,$$
$$a\,\nep^{2r}\leq L< a\,\nep^{8r}\qquad{\rm{if~}}r\geq 1\,.$$
\end{defi}
We let $\mathcal R$ denote the family of all \CZ sets, and observe that  $\mathcal R$ is invariant under left translation. Given  $R \in \mathcal R$, we define its dilated set as $R^*=\{x\in G:~d(x,R)<r\}$. There exists a constant $C_0$ such that $\rho(R^*)\leq C_0\,\rho(R)$ and $R\subset B\big((b_1,b_2,a), C_0r   \big)$. 

 In \cite{HS} it is proved that every integrable function on $G$ admits a \CZ decomposition involving the family $\mathcal R$, and that a new \CZ theory can be developed in this context. By using the \CZ sets, it is natural to introduce an atomic Hardy space $H^1$ on the group $G$, as follows (see \cite{V} for details).

We define an {\emph{atom}} as a function $a$ in $L^1$ such that
\begin{itemize}
\item [(i)] $a$ is supported in a \CZ set $R$;
\item [(ii)]$\|a\|_{\infty}\leq \rho(R)^{-1}\,;$ 
\item [(iii)]$\int a\dir =0$\,.
\end{itemize}
The atomic Hardy space is now defined in a standard way. 
\begin{defi}
The Hardy space $H^{1}$ is the space of all functions $f$ in $ L^1$ which can be written as $f=\sum_j \lambda_j\, a_j$, where $a_j$ are atoms and $\lambda _j$ are complex numbers such that $\sum _j |\lambda _j|<\infty$. We denote by $\|f\|_{H^{1}}$ the infimum of $\sum_j|\lambda_j|$ over such decompositions.
\end{defi}


The new \CZ theory introduced in \cite{HS} is used to study the boundedness of some singular integral operators related to a distinguished Laplacian on $G$, which is defined as follows.

Let $X_0,\,X_1,\,X_2$ denote the left-invariant vector fields  
$$X_0=a\,\partial_a\qquad X_1=a\,\partial_{x_1}\qquad  X_2=a\,\partial_{x_2}\,,$$
which span the Lie algebra $\mathfrak{g}$ of $G$. The Laplacian $\Delta=-(X_0^2+X_1^2+X_2^2)$ is a left-invariant operator which is essentially selfadjoint on $L^2(\rho)$. Since $\Delta$ is positive definite and one-to-one \cite{GQS}, its powers $\Delta^{\alpha}$, $\alpha\in\RR$, have dense domains and are selfadjoint. This makes it possible to form the Riesz transforms of the first order associated with $\Delta$, defined by
\begin{equation}\label{first}
R_i=X_i\,\Delta^{-1/2}\qquad{\rm {and}} \qquad S_i=\Delta^{-1/2}\,X_i,\qquad i=0,1,2\,,
\end{equation}
and the Riesz transforms of the second order, defined by
\begin{equation}\label{second}
R_{ij}=X_iX_j\,\Delta^{-1}\qquad{\rm {and}} \qquad S_{ij}=\Delta^{-1}\,X_iX_j \qquad{\rm{and}} \qquad T_{ij}=X_i\Delta^{-1}X_j\,,
\end{equation}
for $i,j=0,1,2$. The boundedness properties of the Riesz transforms associated with the distinguished Laplacian $\Delta$ defined above have been considered by many authors. Actually, some results in the literature concern the Riesz transforms associated with a distinguished right-invariant Laplacian $\Delta^r$, which is naturally related to $\Delta$ as follows. Let $X_i^r$, $i=0,1,2$, be the right-invariant vector fields on $G$ which agree with $X_i$ at the identity, i.e.,
$$X_0^r=x_1\,\partial_{x_1}+x_2\,\partial_{x_2}+a\,\partial_a\qquad X_1^r=\partial_{x_1}\qquad  X_2^r=\partial_{x_2}\,.$$
It is well known that $X_i^rf=(X_i\check{f})^{\lor}$ for any $f\in C^{\infty}(G)$, where $\check{f}(x)=f(x^{-1})$ for $x\in G$. The Laplacian $\Delta^r=-(X_0^r)^{2}-(X_1^r)^{2}-(X_2^r)^{2}$ is a right-invariant operator which is essentially selfadjoint on $L^2(\lambda)$. We denote by $R_i^r,\,S_i^r, R_{ij}^r,\,S_{ij}^r,\,T_{ij}^r$ the Riesz transforms defined as above by using the right-invariant vector fields and the right-invariant Laplacian instead of the left-invariant ones. It is easy to see that for any $f\in C^{\infty}_c(G)$, $\Delta^rf=(\Delta \check{f})^{\lor}$, 
$$R_i^rf=(R_i \check{f})^{\lor}\,,\qquad S_i^rf=(S_i \check{f})^{\lor}$$
and
$$R_{ij}^rf=(R_{ij} \check{f})^{\lor}\,,\qquad S_{ij}^rf=(S_{ij} \check{f})^{\lor}\,,\qquad T_{ij}^rf=(T_{ij} \check{f})^{\lor}\,.$$
Since $f\to \check{f}$ is an isometry between $L^p(\lambda)$ and $L^p(\rho)$ for $p\in[1,\infty]$, results concerning the boundedness of the right-invariant Riesz transforms with respect to the left Haar measure $\lambda$ may be reformulated in terms of the left-invariant Riesz transforms with respect to the right Haar measure $\rho$. We now summarize some results formulated in terms of the left-invariant Riesz transforms defined by (\ref{first}) and (\ref{second}).

In \cite{GS2, S} G.~Gaudry and P.~Sj\"ogren studied Riesz transforms of the type $X\Delta^{-1/2}$ and $\Delta^{-1/2}X$, where $\Delta$ is a distinguished Laplacian and $X$ is a distinguished vector field, in the context of the group $\RR\ltimes \RR^+$, also known as affine group of the real line. They proved that these operators are of weak type $1$ and bounded on $L^p$, for $1<p<\infty$. In the sequel we sometimes refer to their papers: even if their setting is different, their arguments may be applied also to our context with some slight changes, and so their results carry over.     

Hebisch and Steger then proved that all the operators $R_i$ are of weak type $1$ and bounded on $L^p$ when $1<p\leq 2$ \cite[Theorem 6.4]{HS}. This result was obtained as an application of the \CZ theory on the group $G$. 

The operators $S_i$ are bounded on $L^2$, for $i=0,1,2$. For $i\neq 0$ the operators $S_i$ are of weak type $1$ and bounded on $L^p$ when $1<p\leq 2$, while the operator $S_0$ is not of weak type $1$ but bounded on $L^p$ for $1<p\leq 2$ (Hebisch, private communication).

Since $R_i$ and $S_i$ are bounded on $L^p$ for $p<2$, it follows by duality that $R_i$ and $S_i$ are also bounded on $L^p$ when $2<p<\infty$. 

The second-order Riesz transforms have been studied first in \cite{GQS} in the context of the affine group of the real line, then in \cite{GS1} in the general setting of $NA$ groups of rank $1$, including the group $G$. The operators $T_{ij}$ are of weak type $1$ and bounded on $L^p$ when $1<p<\infty$. The operators $R_{ij}$ and $S_{ij}$ are not of weak type $p$, for any $1\leq p<\infty$.

In this paper we study the $H^1-L^1$ boundedness of the Riesz transforms on the group $G$. Our main results are the following:
\begin{itemize}
\item[(1)] the operators $R_i$, $i=0,1,2$, are bounded from $H^1$ to $L^1$ (Section \ref{Ri}); 
\item[(2)] the operators $S_i$, $i=0,1,2$, are not bounded from $H^1$ to $L^1$ (Sections \ref{Si}, \ref{Szero});
\item[(3)] the operators $T_{ij}$ are bounded from $H^1$ to $L^1$ (Section \ref{Tij}); 
\item[(4)] the operators $S_{ij}$ and $R_{ij}$ are not bounded from $H^1$ to $L^1$ (Sections \ref{Sij}, \ref{Rij}).
\end{itemize}
We remark that since the interpolation spaces between $H^1$ and $L^2$ for the real interpolation method  are the $L^p$ spaces for $1<p<2$ (see \cite{V}), the boundedness of $R_i$ and $T_{ij}$ from $H^1$ to $L^1$ implies their boundedness on $L^p$, for $1<p<2$.

\smallskip
The Riesz operators, and in particular their boundedness on $L^p$ and on the Hardy space $H^1$, have been studied on various Lie groups and Riemannian manifolds. Many results in the literature concern ``doubling spaces'', i.e., measured metric spaces where the volume of balls satisfies the doubling condition. In this context, the Hardy space $H^1$ is defined as in \cite{CW}. 

In the classical setting of $\RR^n$, the Riesz transforms are bounded on $L^p$ for $1<p<\infty$, of weak type $1$ and bounded on $H^1$ \cite[III.3]{St}.

For nilpotent Lie groups and first-order Riesz operators, the $L^p$-bounded\-ness, for $1<p<\infty$, the weak type $1$ and the $H^1$-boundedness were proved by N.~Lohou\'e and N.~Varopoulos \cite{LV}. Subsequently, this was extended  to all connected Lie groups of polynomial growth by L.~Saloff-Coste \cite{SC} and G.~Alexopoulos \cite{A}. 

In the setting of symmetric spaces of noncompact type, J.-P.~Anker \cite{A2} considered Riesz transforms associated with the Laplace--Beltrami operator. He proved the weak type $1$ estimate for the first-order operators and the $L^p$-estimates for operators of arbitrary order.

On a Riemannian manifold the Riesz transform $R=\nabla \Delta^{-1/2}$, where $\nabla$ is the gradient and $\Delta$ is the Laplace--Beltrami operator, has been considered. If the manifold has nonnegative Ricci curvature, then the Riesz transform $R$ is bounded on $L^p$, $1<p<\infty$, of weak type $1$ and bounded from $H^1$ to $L^1$ \cite{B, CL}. Subsequently, T.~Coulhon and X.T.~Duong proved that on a Riemannian manifold with the doubling property whose heat kernel verifies an upper estimate on the diagonal, $R$ is of weak type $1$ and bounded on $L^p$, for $1<p\leq 2$ \cite{CD1}. The connection between the $L^p$-boundedness of the Riesz transform, Poincar\'e inequalities and heat kernel estimates is also studied in \cite{AC, ACDH, CD2, CL2}. In Riemannian manifolds satisfying the doubling condition and the Poincar\'e inequality, E. Russ \cite{R} proved that $R$ is bounded from $H^1$ to $L^1$; then M. Marias and Russ \cite{MR} proved the boundedness on $H^1$ of the linearized 
 Riesz transforms.

The previous results do not apply to our space $(G,d,\rho)$, since it is of exponential growth and the doubling condition fails.

\vspace{.5cm}

Our paper is organized as follows: Section \ref{kernels} contains the analysis of the kernels of the Riesz transforms. The $H^1-L^1$-boundedness of the operators $R_i$ is proved in Section \ref{Ri}, as a consequence of a general boundedness theorem for integral operators. In Section~\ref{Si}, we prove the unboundedness from $H^1$ to $L^1$ of the operators $S_1$ and $S_2$, and in Section \ref{Szero} that of $S_0$. We analyze the local part of the second-order Riesz transforms in Section \ref{local}, proving that they are bounded from $H^1$ to $L^1$. In Section \ref{Tij} we show that the operators $T_{ij}$ are bounded from $H^1$ to $L^1$. Finally, we show that the global part of the operators $S_{ij}$ and $R_{ij}$ are not bounded from $H^1$ to $L^1$ in Sections \ref{Sij} and \ref {Rij}.

\smallskip
In the following, $C$ denotes a positive, finite constant which may vary from line to line and may depend on parameters according to the context. Given two quantities $f$ and $g$, by $f\sim g$ we mean that there exists a constant $C$ such that $1/C\leq f/g\leq C$.


\section{The convolution kernels of the Riesz transforms}\label{kernels}
In this section, we analyze the convolution kernels of the Riesz transforms of the first and the second order. First recall that the definition of the convolution of two functions $f,\,g$ on $G$ is 
\begin{align*}
f\ast g(x)&=\int_Gf(xy^{-1})\,g(y)\dir (y)\qquad\forall x\in G\,.
\end{align*}
Let $V$ denote the space $\{\Delta u:~u\in C^{\infty}_c(G)\}$. In \cite{GS1} it is verified that $V$ is a dense subspace of $L^2$ and that $V\subset D(\Delta^{-1})\subset D(\Delta^{-1/2})$. For $\alpha>0$, we denote by $U_{\alpha}$ the convolution kernel of $\Delta^{-\alpha/2}$, in the sense that $\Delta^{-\alpha/2}f=f\ast U_{\alpha}$, for all $f\in V$. Since
$$\Delta^{-\alpha/2}=\frac{1}{\Gamma(\alpha/2) }\int_0^{\infty}t^{\alpha/2-1}\nep^{-t\Delta}\di t\,,$$
we have that
$$U_{\alpha}=\frac{1}{\Gamma(\alpha/2) }\int_0^{\infty}t^{\alpha/2-1}p_t\di t\,,$$
where $p_t$ denotes the heat kernel of $\Delta$. It is well known \cite[Theorem 5.3, Proposition 5.4]{CGGM}, \cite[Formula (5.7)]{ADY} that 
$$p_t(x)=\frac{1}{8\pi^{3/2}}\,\delta^{1/2}(x)\,\frac{r(x)}{\sinh r(x)}\,t^{-3/2}\,\nep^{-\frac{r^2(x)}{4t}}   \qquad \forall x\in G\,,$$
where $r(x)$ denotes as before the distance of $x$ from the identity. Hence, for $\alpha<3$
\begin{align*}
U_{\alpha}(x)&=\frac{1}{\Gamma(\alpha/2) }\,\frac{1}{8\pi^{3/2}}\,\delta^{1/2}(x)\,\frac{r(x)}{\sinh r(x)}\,\int_0^{\infty}t^{\alpha/2-1}t^{-3/2}\,\nep^{-\frac{r^2(x)}{4t}} \di t\\
&=\frac{1}{\Gamma(\alpha/2) }\,\frac{2^{1-\alpha}}{\pi^{3/2}}\,\delta^{1/2}(x)\,\frac{r(x)}{\sinh r(x)}\int_0^{\infty}r(x)^{\alpha-3}v^{2-\alpha}\nep^{-v^2}\di v\\
&=C_{\alpha}\,\delta^{1/2}(x)\,\frac{r^{\alpha-2}(x)}{\sinh r(x)}\qquad \forall x\in G\,.
\end{align*}
We consider the cases $\alpha=1$ and $\alpha=2$ and get that $C_1=\frac{1}{2\pi^2}$ and $C_2=\frac{1}{4\pi}$. We denote by $U=U_1$ the convolution kernel of $\Delta^{-1/2}$ given by
\begin{equation}\label{formulaU}
U(x)=\frac{1}{2\pi^{2}}\,\delta^{1/2}(x)\,\frac{1}{r(x)\,\sinh r(x)}\qquad \forall x\in G,
\end{equation}
and by $W=U_2$ the convolution kernel of $\Delta^{-1}$ given by
\begin{equation}\label{formulaW}
W(x)=\frac{1}{4\pi}\,\delta^{1/2}(x)\,\frac{1}{\sinh r(x)}\qquad \forall x\in G\,.
\end{equation}
Since $R_i=X_i\,\Delta^{-1/2}$, we get for all $f\in V$ and $x\in G$
\begin{align}\label{pvXiU}
R_if(x)&=X_i(f\ast U)(x)=\int X_{i,x}f(xy^{-1})\,U(y)\dir(y)\nonumber\\
&=\lim_{\varepsilon\to 0}\int_{r(y)>\varepsilon } X_{i,x}f(xy^{-1})\,U(y)\dir(y)\nonumber\\
&=-\lim_{\varepsilon\to 0}\int_{r(y)>\varepsilon } X_{i,y}f(xy^{-1})\,U(y)\dir(y)\nonumber\\
&=\lim_{\varepsilon\to 0}\int_{r(y)>\varepsilon } f(xy^{-1})\,X_{i,y}U(y)\dir(y)\,,
\end{align}
where the last step follows by integration by parts, as in \cite[Section 3]{S}. Thus the convolution kernel of $R_i$ is the distribution ${\rm{pv~}} k_i$, where $k_i=X_iU$. Moreover, for $f\in V$ and $x\notin {\rm{supp}} f$
\begin{align}\label{Ri(x,y)}
R_if(x)&=\int_G f(xy)\,k_i(y^{-1})\dil(y)\nonumber\\
&=\int_G f(y)\,k_i(y^{-1}x)\,\delta(y)\dir(y)\nonumber\\
&=\int_G f(y)\,R_i(x,y)\dir(y)\,,
\end{align}
where $R_i(\cdot,\cdot)$ denotes the integral kernel of $R_i$, related to $k_i$ by
\begin{equation}\label{kernelRi}
R_i(x,y)=\delta(y)\,k_i(y^{-1}x)\qquad \forall x,\,y\in G,\qquad x\neq y\,.
\end{equation}
We now consider the operators $S_i$. By arguing as in \cite[page 246-247]{GS2}, it is easy to see that if $f\in C^{\infty}_c(G)$, then $X_if\in D(\Delta^{-1/2})$, so that $S_i$ is well defined on $C^{\infty}_c(G)$. Moreover, for all $f\in C^{\infty}_c(G)$ and $g\in V$
$$\langle  S_if,g\rangle =\langle  \Delta^{-1/2}\,X_if,g\rangle=\langle X_if,\Delta^{-1/2}g\rangle=-\langle f, X_i\,\Delta^{-1/2}g\rangle=-\langle f, R_ig\rangle\,.$$
Thus by (\ref{kernelRi}) we deduce that the integral kernel of $S_i$ is given by
\begin{equation}\label{kernelSi}
S_i(x,y)=-\overline{R_i(y,x)}=-\delta(x)\,k_i(x^{-1}y)\qquad\forall x,\,y\in G,\qquad x\neq y\,.
\end{equation}
We now compute $k_i$ explicitly. To do so, we shall need the following simple lemma.
\begin{lemm}\label{derivativedistance}
At any point $(x_1,x_2,a)\neq (0,0,1)$ in $G$, the derivatives of $r$ along the vector fields $X_i$ are given by
$$X_ir(x_1,x_2,a)=\begin{cases}
\frac{a-a^{-1}-a^{-1}(x_1^2+x_2^2)}{2\,\sinh r(x_1,x_2,a)}=\frac{a-\cosh r}{\sinh r}&{\rm{if}}~i=0\\
\frac{x_i}{\sinh r(x_1,x_2,a)}&{\rm{if}}~i=1,2\,.
\end{cases}$$
\end{lemm}
\begin{proof}
It suffices to differentiate the expression
$$\cosh r(x_1,x_2,a)=\frac{a+a^{-1}+a^{-1}(x_1^2+x_2^2)}{2}\,,$$
with respect to $X_i$. For $X_0=a\,\partial_a$ we obtain 
$$\sinh r(x_1,x_2,a)\,X_0r(x_1,x_2,a)=a~\frac{1-a^{-2}-a^{-2}(x_1^2+x_2^2)}{2}\,,$$
which gives the result for $i=0$. The cases of $X_i=a\,\partial_{x_i}$, $i=1,2$, are similar. 
\end{proof}
By (\ref{formulaU}) and Lemma \ref{derivativedistance} for $i=1,2$ and $(x_1,x_2,a)\neq (0,0,1)$ we get
\begin{align}\label{k1}
k_i(x_1,x_2,a)&=X_iU(x_1,x_2,a)\nonumber\\
&=-\frac{1}{2\,\pi^{2}}\,a^{-1}\,\frac{\sinh r+r\,\cosh r}{r^2\,\sinh^2 r}\,X_ir(x_1,x_2,a)\nonumber\\
&= -\frac{1}{2\,\pi^{2}}\,a^{-1}\,x_i\,\frac{\sinh r+r\,\cosh r}{r^2\,\sinh^3 r}\,.
\end{align}
For $i=0$ and $(x_1,x_2,a)\neq (0,0,1)$ we get
\begin{align}\label{k0}
&k_0(x_1,x_2,a)=X_0U(x_1,x_2,a)\nonumber\\
&=\frac{1}{2\,\pi^{2}}\Big[-\,a\,a^{-2}\,\frac{1}{r\,\sinh r}-\,a^{-1}\,\frac{\sinh r+r\,\cosh r}{r^2\,\sinh^2 r}\,X_0r(x_1,x_2,a)\Big]\nonumber\\
&=\frac{1}{2\,\pi^{2}}\Big[ -\,a^{-1}\,\frac{1}{r\,\sinh r}-\,a^{-1}\,\frac{a-a^{-1}-a^{-1}(x_1^2+x_2^2)}{2}\,\frac{\sinh r+r\,\cosh r}{r^2\,\sinh^3 r}\Big]\nonumber\\
&= -U(x_1,x_2,a) +\frac{1}{2\,\pi^{2}}\,\,\frac{-1+a^{-2}+a^{-2}(x_1^2+x_2^2)}{2}\,\frac{\sinh r+r\,\cosh r}{r^2\,\sinh^3 r}\,.
\end{align}

\smallskip
We now consider the second-order Riesz transforms. We shall regard $\Delta^{-1}$ as the operator of convolution by the kernel $W$. The operators $R_{ij}=X_iX_j\Delta^{-1}, \,S_{ij}=\Delta^{-1}X_iX_j,\,T_{ij}=X_i\Delta^{-1}X_j$ are then properly defined on $C^{\infty}_c(G)$, with values in $C^{\infty}(G)$. By arguing as in \cite[Lemma 6]{GS1} we may show that there exist distributions $k_{ij},\, \ell_{ij}, \,g_{ij}$ such that for any $f\in C^{\infty}_c(G)$
$$R_{ij}f=f\ast k_{ij}\qquad S_{ij}f=f\ast \ell_{ij}\qquad   T_{ij}f=f\ast g_{ij}\,.$$
To compute these convolution kernels, we recall some simple properties of right- and left-invariant vector fields, which are the analogs of those proved in \cite[Section 4.2]{GS1} with respect to the measure $\lambda$.

Given a vector $Z$ in $\mathfrak{g}$, we here denote by $Z^r$ and $Z^{\ell}$ the right-invariant and left-invariant vector fields on $G$ which agree with $Z$ at the identity, defined by
$$Z^rf(x)=\frac{d}{dt}\Big\lvert_{ t=0}f\big({\rm{exp}}(tZ)\,x  \big)\qquad{\rm{and}}\qquad Z^{\ell}f(x)=\frac{d}{dt}\Big\lvert_{ t=0}f\big(x\,{\rm{exp}}(tZ) \big)$$
for any $f\in C^{\infty}(G)$ and $x\in G$. It is easy to check that for any $f\in C^{\infty}(G)$
\begin{equation}\label{ZrZell}
Z^r\check{f}=(Z^{\ell}f)^{\lor}\,.
\end{equation}
Let $k$ be  a distribution on $G$, $f\in C_c(G)$ and $g\in C^{\infty}_c(G)$. Then
\begin{equation}\label{astscalar}
\langle f\ast k,g\rangle= \langle k,\check{f}\ast g \rangle\,.
\end{equation}
The left-invariant derivative $Z^{\ell}k$ of $k$ is the distribution such that for any $g\in C^{\infty}_c(G)$
\begin{equation}\label{meno}
\langle Z^{\ell} k,g\rangle=- \langle k,Z^{\ell}g \rangle  \,.
\end{equation}
If $f,\,g \in C^{\infty}_c(G)$, then
$$\langle Z^rf,g\rangle= \langle f,-Z^rg \rangle+Z\delta(e)\,\langle f,g \rangle  \,.$$
So it is natural to define the right-invariant derivative of a distribution $k$ as the distribution $Z^rk$ for which
\begin{equation}\label{Zrkg}
\langle Z^rk,g\rangle=\langle k,-Z^rg \rangle+Z\delta(e)\,\langle k,g \rangle  \qquad \forall g\in C^{\infty}_c(G)\,.
\end{equation}
It is easy to verify that
\begin{equation}\label{spostacampi}
Z^{\ell}\big(f\ast k\big)=f\ast Z^{\ell}k\qquad{\rm{and}}\qquad Z^r\big(f\ast k\big)=Z^rf\ast k\,.
\end{equation}
By (\ref{astscalar}) and (\ref{Zrkg}) we deduce that
\begin{equation}\label{ZrZl}
Z^{\ell}f\ast k=f\ast \big(-Z^rk+Z\delta(e)k  \big)\,.
\end{equation}
Applying (\ref{spostacampi}) we get that for any $f\in C^{\infty}_c(G)$ 
$$R_{ij}f=X_iX_j\Delta^{-1}f=X_iX_j(f\ast W)=X_i(f\ast X_jW)=f\ast X_iX_jW\,.$$
Thus the convolution kernel of $R_{ij}$ is 
\begin{equation}\label{formulakij}
k_{ij}=X_iX_jW\,,
\end{equation}
the derivative taken in the distribution sense. We denote by $R_{ij}(\cdot,\cdot)$ the integral kernel of $R_{ij}$ defined by $R_{ij}(x,y)=\delta(y)\,k_{ij}(y^{-1}x)$, for $x\neq y$.

Moreover, by (\ref{meno}) for all $f,\,g\in C^{\infty}_c(G)$
\begin{align*}
\langle  S_{ij}f,g\rangle &=\langle  \Delta^{-1}\,X_iX_jf,g\rangle=\langle X_iX_jf,\Delta^{-1}g\rangle\\
&=-\langle X_jf, X_i\,\Delta^{-1/2}g\rangle=\langle f, R_{ji}g\rangle\,.
\end{align*}
This implies that the integral kernel of $S_{ij}$ is given by
\begin{equation}\label{kernelSij}
S_{ij}(x,y)=R_{ji}(y,x)=\delta(x)\,k_{ji}(x^{-1}y)  \qquad \forall  x\neq y\,.
\end{equation}
It easily follows that the convolution kernel of $S_{ij}$ is $\ell_{ij}=\delta\check{k}_{ji}$. 

Applying (\ref{ZrZl}), we get that for any $f\in C^{\infty}_c(G)$
\begin{align*}
T_{ij}f&=X_i\Delta^{-1}X_jf=X_i\big(X_jf\ast W\big)\\
&=X_jf\ast X_iW=f\ast  \big(-X_j^rX_iW+X_j\delta(e)X_iW  \big)\,.
\end{align*}
Thus the convolution kernel of $T_{ij}$ is $g_{ij}=-X_j^rX_iW+X_j\delta(e)X_iW$.

To avoid long computations, we do not compute explicitly the kernels of the second-order Riesz transforms, but we shall find their behavior away from the identity, i.e., in the complement of the unit ball $B_1$. 

In the sequel, we shall denote by $R(r)$ any series of the type $\sum_{k=1}^{\infty}c_k\nep^{-2kr}$, where the $c_k$ are real numbers and the series is convergent for $r> 1$; we denote by $S(r)$ any function of the type $1+R(r)$. These functions may vary from occurrence to occurrence. It is easy to see that a function $S(r)=1+R(r)$ may be differentiated termwise and its derivative is $S'(r)=R(r)$. Moreover, multiplying two functions of the type $S$, we obtain a function of the same kind. 

Lemma \ref{derivativedistance} implies that for points $(x_1,x_2,a)\in \overline{B}_1^c$
\begin{equation}\label{Xir}
X_ir(x_1,x_2,a)=\begin{cases}
2x_i\,\nep^{-r}\,S(r)&{\rm{if~}}i=1,2\\
2a\,\nep^{-r}\,S(r)-S(r)&{\rm{if~}}i=0\,,
\end{cases}
\end{equation}
and by (\ref{formulaW})
\begin{equation}\label{Wexpansion}
W(x_1,x_2,a)=\frac{1}{2\pi}\,a^{-1}\,\nep^{-r}\,S(r)\,.
\end{equation}
Let $\ZZ^3_+$ be the set of $m=(m_0,m_1,m_2)$ in $\ZZ^3$ such that $m_1,m_2\geq 0$ and $m_0\geq -1$. We denote by $|m|$ the sum $m_0+m_1+m_2$ and by $x^m$ the product $x_1^{m_1}x_2^{m_2}a^{m_0}$. The principal term of $W$ is of the type $x^{m}\nep^{-p r}$, where $|m|-p=-2$. We shall study the integrability of similar expressions in an elementary lemma, and first split $\overline{B}_1^c$ into two parts, as follows:
$$G_+=\{(x_1,x_2,a)\in \overline{B}_1^c:~a>1\}\,,$$
and
$$G_-=\{(x_1,x_2,a)\in \overline{B}_1^c:~a<1\}\,.$$
\begin{lemm}\label{integrability and derivative}
Let $m$ be in $\ZZ^3_+$ and $p\in \NN$. 
\begin{itemize}
\item[(i)] The function $ x^{m}\,\nep^{-p\, r}$ is integrable in $G_+$ if and only if
\begin{equation}\label{conditionsa>1}
m_1+m_2-2p< -2\qquad {\rm{and}}\qquad |m|-p<  -2  \,.
\end{equation}
\item[(ii)] The function $ x^{m}\,\nep^{-p\, r}$ is integrable in $G_-$ if and only if
\begin{equation}\label{conditionsa<1}
m_1+m_2-2p< -2\qquad {\rm{and}}\qquad m_0+p> 0\,.
\end{equation}
\end{itemize}
\end{lemm}
\begin{proof}
If $(x_1,x_2,a)\in G_+$, then $\nep^{r} \sim a\,(1+a^{-2}\,|(x_1,x_2)|^2)$, so that
\begin{align*}
&\int_{G_+}  a^{m_0}\,|x_1|^{m_1}\,|x_2|^{m_2}\,\nep^{-p\, r}\di x_1\di x_2 \frac{\di a}{a}\\
&\sim   \int_1^{\infty} a^{m_0}\int_{\RR^2} |x_1|^{m_1}\,|x_2|^{m_2}\,a^{-p}\,(1+a^{-2}\,|(x_1,x_2)|^2)^{-p}\di x_1\di x_2 \frac{\di a}{a}\,.
\end{align*}
Under the change of variables $a^{-1}(x_1,x_2)=(y_1,y_2)$, this transforms into the product
\begin{align*}
&\int_1^{\infty} a^{|m|-p+1}\di a \; \int_{\RR^2} |y_1|^{m_1}\,|y_2|^{m_2}\,(1+|(y_1,y_2)|^2)^{-p}\di y_1\di y_2\,.
\end{align*}
Here the integral in $a$ converges if and only if $|m|-p< -2$. By means of polar coordinates, the second integral is seen to converge if and only if $m_1+m_2-2p< -2$. This proves (i).

 If $(x_1,x_2,a)\in G_-$, then $\nep^{r}\sim  a^{-1}\,(1+|(x_1,x_2)|^2)$, so that
\begin{align*}
&\int_{G_-}  a^{m_0}\,|x_1|^{m_1}\,|x_2|^{m_2}\,\nep^{-p\, r}\di x_1\di x_2 \frac{\di a}{a}\\
&\sim  \int_0^{1} a^{m_0+p-1} {\di a} \; \int_{\RR^2} |x_1|^{m_1}\,|x_2|^{m_2}\,(1+|(x_1,x_2)|^2)^{-p}\di x_1\di x_2\,,
\end{align*}
and (ii) follows.
\end{proof}
To study the higher order derivatives of $W$, we start with the derivatives along $X_0,X_1,X_2$ of an expression $x^m\nep^{-pr}\,S(r)$, as above. We shall always have
\begin{equation}\label{conditions2}
m_1+m_2-2p\leq -2\qquad |m|-p \leq  -2 \qquad {\rm{and}}\qquad m_0+p\geq  0\,,
\end{equation}
which does not imply the integrability of $x^m\nep^{-pr}$. For many remainder terms, we shall denote by $Q(x)$ any finite sum of terms $x^n\nep^{-qr}\,R(r)$, where $|n|-q\leq -2$, $n_0+q\geq 0$ and $n_1+n_2-2q\leq  -2$, so that $Q(x)$ is integrable in $B_1^c$. By (\ref{Xir}) we get that in $\overline{B}_1^c$
\begin{align}\label{X1term}
X_1\big( x^m\,\nep^{-pr}\,S(r) \big)&= m_1\,a^{m_0+1}x_1^{m_1-1}x_2^{m_2}\nep^{-pr}\,S(r)-\\
& -p\, a^{m_0}x_1^{m_1}x_2^{m_2}\nep^{-pr}\,2x_1\,\nep^{-r}\,S(r)+\nonumber\\
&+a^{m_0}x_1^{m_1}x_2^{m_2}\nep^{-pr}\,\,R(r)\,2x_1\,\nep^{-r}\,S(r)\nonumber\\
&=m_1\,a^{m_0+1}x_1^{m_1-1}x_2^{m_2}\nep^{-pr} -2p\, a^{m_0}x_1^{m_1+1}x_2^{m_2}\nep^{-(p+1)r}+\nonumber\\
&+Q(x)\,.
\end{align}
By symmetry an analogous formula holds for $i=2$. From (\ref{Xir}) we get
\begin{align}\label{X0term}
X_0\big(x^m\,\nep^{-pr} \,S(r)\big)&= m_0 a^{m_0}x_1^{m_1}x_2^{m_2}\,\nep^{-pr}\,S(r)-\\
& -p\,a^{m_0}x_1^{m_1}x_2^{m_2}\,\nep^{-pr}\,\big[2a\,\nep^{-r}\,S(r)-S(r)  \big]+\nonumber\\
&+ a^{m_0}x_1^{m_1}x_2^{m_2}\,\nep^{-pr}\,R(r)\big[2a\,\nep^{-r}\,S(r)-S(r)  \big] \nonumber  \\
&=(m_0+p)\,  a^{m_0}x_1^{m_1}x_2^{m_2}    \,\nep^{-pr}- 2p\,a^{m_0+1}x_1^{m_1}x_2^{m_2}\nep^{-(p+1)r}+\nonumber\\
&+Q(x)\,.    
\end{align}
Differentiating the expression (\ref{Wexpansion}) for $W$ and applying (\ref{X1term}) and (\ref{X0term}), we get that in $\overline{B}_1^c$
\begin{align*}
X_jW(x_1,x_2,a)&=-\frac{1}{2\pi }\,a^{-1}x_j\,\nep^{-2r}+Q(x)\qquad {\rm{if~}}j=1,2,
\end{align*}
and
\begin{align*}
X_0W(x_1,x_2,a)&=-\frac{1}{2\pi}\,\nep^{-2r}+Q(x)\,.
\end{align*}
We now differentiate $W$ a second time, applying again (\ref{X1term}) and (\ref{X0term}) and also the observation that $X_jQ(x)=Q(x)$, for $i,j=0,1,2$. The result is that there exist constants $\alpha_{ij},\,\beta_{ij}\in\RR$ and $m,\,n\in \ZZ^3_+$ such that in $\overline{B}_1^c$
\begin{align}\label{kappaij}
k_{ij}(x)&=\alpha_{ij}\,x^m\,\nep^{-2r} +\beta_{ij}\,x^n\,\nep^{-3r} +Q(x)\,,
\end{align}
where $\beta_{ij}\neq 0$, $|m|=0,\,|n|=1$, $m_1+m_2-4< -2$, $n_1+n_2-6< -2$, $m_0+2>0$ and $n_0+3>0$. This means that $k_{ij}$ has a principal part in $\overline{B}_1^c$ given by at most two nonintegrable terms, while the remaining part of the kernel is integrable. Finally, we estimate the derivative of $k_{ij}$ along the vector field $X_2$. We get that, for $i,j=0,1,2$, there exist constants $\gamma_{ij},\,\eta_{ij},\,\sigma_{ij},\,\theta_{ij}\in\RR$ and $h,\,\ell,\,m,\,n\in \ZZ^3_+$ such that in $\overline{B}_1^c$
\begin{align}\label{X2kappaij}
X_2k_{ij}(x)&=\gamma_{ij}\,x^h\,\nep^{-2r} +\eta_{ij}\,x^{\ell}\,\nep^{-3r} +  \sigma_{ij}\,x^{m}\,\nep^{-3r} +\theta_{ij}\,x^{n}\,\nep^{-4r}       +Q(x)\,,
\end{align}
where $\theta_{ij}\neq 0$, $|h|=0,\,|\ell|=|m|=1,\,|n|=2$.

\section{$H^1-L^1$-boundedness of $R_i$}\label{Ri}
In this section we prove that the Riesz transforms $R_i$ are bounded from $H^1$ to $L^1$, for $i=0,1,2$. 

This result is a consequence of the following boundedness theorem for integral operators. Note that the hypotheses of the following proposition are the same as those of \cite[Theorem 2.1]{HS}.
\begin{prop}\label{Teolim}
Let $T$ be a linear operator bounded on $L^2$ such that $T=\sum_{j\in \ZZ}T_j$, where
\begin{itemize}
\item [(i)] the series converges in the strong operator topology of $L^2$;
\item [(ii)] every $T_j$ is an integral operator with integral kernel $T_j$;
\item [(iii)] there exist positive constants $\alpha ,A,\varepsilon$ and $c>1$ such that
\begin{align}\label{stime1}
\int_G|T_j(x,y)|\,\big(1+c^j d(x,y)\big)^{\varepsilon}\,{\rm{d}}\rho (x) &\leq A\qquad\forall y\in G;
\end{align}
\begin{align}\label{stime2}
\int_G|T_j(x,y)-T_j(x,z)| \,{\rm{d}}\rho (x) &\leq A\,\big(c^j d(y,z)\big)^{\alpha}\qquad\forall y,z\in G\,.
\end{align}
\end{itemize}
Then $T$ is bounded from $H^1$ to $L^1$.
\end{prop}
\begin{proof}
We first show that there exists a constant $C$ such that for any atom $a$
\begin{equation}\label{Ta}
\|Ta\|_1\leq C\,.
\end{equation}
Let $R$ be the \CZ set supporting $a$, and denote by $c_R$ the center of $R$ and by $R^*$ its dilated set (defined in Section \ref{intro}). We estimate the integral of $Ta$ on $R^*$ by the Cauchy--Schwarz inequality:
\begin{equation}\label{suR^*}
\int_{R^*} |Ta|\dir\leq \|Ta\|_{2}\,\rho(R^*)^{1/2} \leq C\,|\!|\!| T |\!|\!|_{L^2\to L^2}\,\|a\|_2\,\rho(R)^{1/2}\leq C\,|\!|\!| T |\!|\!|_{L^2\to L^2}\,.
\end{equation}
It is easy to show that from the estimates (\ref{stime1}) and (\ref{stime2}) it follows that
\begin{equation}\label{Hormander}
\sup_{R\in\mathcal R}\sup_{y,\,z\in R}\int_{(R^*)^c}|T(x,y)-T(x,z)  |\dir(x)<\infty\,,
\end{equation}
where $T$ is the integral kernel of $T$. Thus the integral of $Ta$ on the complementary set of $R^*$ is estimated as follows: 
\begin{align*}
\int_{(R^{*})^c} |Ta|\dir&\leq \int_{(R^{*})^c}\Big|\int_R T(x,y)\,a(y)\dir(y)  \Big|\dir(x)\\
&=
\int_{(R^{*})^c}\Big|\int_R [T(x,y)-T(x,c_R)]\,a(y)\dir(y)  \Big|\dir(x)\\
&\leq \int_{(R^{*})^c}\int_R |T(x,y)-T(x,c_R)|\,|a(y)|\dir(y)\dir(x)\\
&=\int_R|a(y)|\Big( \int_{(R^{*})^c} |T(x,y)-T(x,c_R)|\dir(x) \Big)\dir(y)\\
&\leq  \|a\|_1\,\sup_{y\in R}\int_{(R^{*})^c}|T(x,y)-T(x,c_R)|\dir(x)\\
&\leq  C\,.
\end{align*}
This concludes the proof of (\ref{Ta}). We shall deduce from (\ref{Ta}) that $T$ is bounded from $H^1$ to $L^1$. Indeed, by \cite[Theorem 2.1]{HS} $T$ is bounded from $L^1$ to $L^{1,\infty}$. Now take a function $f$ in $H^1$ and suppose that $f=\sum_{j=1}^{\infty }\lambda_ja_j$ is an atomic decomposition with $\sum_j|\lambda_j|\sim \|f\|_{H^1}$. Define $f_N=\sum_{j=1}^{N}\lambda_ja_j$. Since $f_N$ converges to $f$ in $L^1$, $Tf_N=\sum_{j=1}^N\lambda_jTa_j$ converges to $Tf$ in $L^{1,\infty}$. On the other hand, by (\ref{Ta})
$$\|Tf_N-\sum_{j=1}^{\infty} \lambda_jTa_j\|_1\leq \sum_{j=N+1}^{\infty}|\lambda_j|\,\|Ta_j\|_1 \leq C\,\sum_{j=N+1}^{\infty}|\lambda_j|\,,$$
so that $Tf_N$ converges to $\sum_{j=1}^{\infty}\lambda_jTa_j$ in $L^1$. This implies that $Tf=\sum_{j=1}^{\infty}\lambda_jTa_j\in L^1$ and $\|Tf\|_1\leq C\,\|f\|_{H^1}\,$.
\end{proof}
We now easily obtain the following theorem.
\begin{theo}
The Riesz transforms $R_i$, for $i=0,1,2$, are bounded from $H^1$ to $L^1$.
\end{theo}
\begin{proof}
In the proof of \cite[Theorem 6.4]{HS}, it is shown that the operator $R_i$ satisfies the assumptions of Proposition \ref{Teolim}. Thus $R_i$ is bounded from $H^1$ to $L^1$.
\end{proof}


\section{Unboundedness of $S_1$ and $S_2$}\label{Si}
In this section we prove that the Riesz transforms $S_1$ and $S_2$ are not bounded from $H^1$ to $L^1$. To do so, we shall define an atom $a$ on $G$ such that the images of $a$ under these operators are not integrable in a region far from the support of the atom (see Theorem \ref{S1}). By symmetry it suffices to consider the case $i=1$.

Differentiating the expression (\ref{k1}) for $k_1$ along the vector field $X_2$ and applying Lemma \ref{derivativedistance}, we obtain that
\begin{align}\label{X2k1}
&X_2k_1(x_1,x_2,a)=-\frac{1}{2\,\pi^{2}}\,a^{-1}\,x_1\,X_2r(x_1,x_2,a)\Big[\frac{r^2\,\sinh^3 r(2\cosh r+r\sinh r)}{r^4\,\sinh^6 r}-\nonumber\\
&-\frac{(\sinh r+r\cosh r)(2r\sinh^3 r+3r^2\,\sinh^2r\cosh r)}{r^4\,\sinh^6 r}\Big]\nonumber\\
&=\frac{1}{2\,\pi^{2}}\,a^{-1}\,~\frac{x_1\,x_2}{\sinh r}~\frac{2\,r^2\,\cosh^2 r+r^2+2\sinh^2 r+3r\,\sinh r\cosh r}{r^3\,\sinh^4 r}\,.
\end{align}
\begin{lemm}\label{estimateX2k1}
There exist regions $\Gamma''\subset  \Gamma'\subset  \Gamma$, a positive continuous function $\Phi$ on ${\Gamma}$ and a positive constant $C$ such that 
\begin{itemize}
\item[(i)] $X_2k_1\geq C\,\Phi$ in ${\Gamma}$;
\item[(ii)] for any $(x_1,x_2,a)$ in $\Gamma'$ and $\tau$ in $[0,1/4]$, the point $(x_1,x_2,a)\cdot(0,\tau,1)$ is in ${\Gamma}$ and
$$\Phi\big((x_1,x_2,a)\cdot(0,\tau,1)\big)= \Phi(x_1,x_2,a) \,;$$
\item[(iii)] $\int_{\Gamma''}\Phi \dir=\infty$\,.
\end{itemize}
Let $E$ be the parallelepiped $(-1/2,1/2)\times (-1/4,0)\times (1,2)$. Then 
\begin{equation}\label{inclusion}
\Gamma''\cdot E^{-1}\cdot E\subseteq \Gamma'\,.
\end{equation}
\end{lemm}
\begin{proof}
Given $B,\,A>1$ and $0<\varepsilon <1$ we define $\Gamma$ as the region 
$$\Gamma=\{(x_1,x_2,a)\in G:~1-\varepsilon <x_2/x_1<1+\varepsilon , \,x_1>Ba,\,a>A\}\,.$$
For any $(x_1,x_2,a)$ in ${\Gamma}$ we have that
$$\frac{a^{-1}x_1^2}{2}<\cosh r(x_1,x_2,a)<C \,a^{-1}\,x_1^2\,,$$
and $r(x_1,x_2,a)>1$, for $A$ and $B$ sufficiently large. Moreover, since $\nep^{r}<2\,\cosh r< C\,a^{-1}\,x_1^2$,  
\begin{align*}
r(x_1,x_2,a)&\leq C\,\log(a^{-1}\,x_1^2)\,.
\end{align*}
By the formula (\ref{X2k1}) it is clear that $X_2k_1$ is positive on $\Gamma$. Considering the first term in the numerator of the last fraction in (\ref{X2k1}), we see that for $(x_1,x_2,a)$ in $\Gamma$
\begin{align*}
X_2k_1(x_1,x_2,a)&\geq C\,a^{-1}\,x_1^2\,\frac{1}{r\,\cosh^3 r}\geq C~\frac{a^{-1}\,x_1^2}{\log(a^{-1}\,x_1^2)\,(a^{-1}\,x_1^2)^3}\,.
\end{align*}
We define 
$$\Phi(x_1,x_2,a)=\frac{1}{\log(a^{-1}\,x_1^2)\,(a^{-1}\,x_1^2)^2}\,.$$
The condition (i) is verified. We now define
\begin{align*} 
\Gamma'&=\{(x_1,x_2,a)\in \Gamma:~1-\varepsilon' <x_2/x_1<1+\varepsilon' , \,x_1>B'a\}\,,\\
\Gamma''&=\{(x_1,x_2,a)\in \Gamma:~1-\varepsilon'' <x_2/x_1<1+\varepsilon'' , \,x_1>B''a,\,a>2A\}\,,
\end{align*}
where $B''>B'>B,\,0<\varepsilon''<\varepsilon'<\varepsilon<1$ have to be chosen. 

Let $(x_1,x_2,a)$ be a point in $\Gamma'$ and $\tau$ in $[0,1/4]$. Then $(x_1,x_2,a)\cdot(0,\tau,1)=(x_1, x_2+a\,\tau,a)\,.$ It is easy to see that we may choose $B',\varepsilon'$ such that $(x_1, x_2+a\,\tau,a)\in \Gamma$. Moreover,
\begin{align*}
\Phi\big((x_1,x_2,a)\cdot(0,\tau,1)\big)&= \frac{1}{\log(a^{-1}\,x_1^2)\,(a^{-1}\,x_1^2)^2} = \Phi(x_1,x_2,a)\,,
\end{align*}
as required in (ii). To prove (iii), we integrate $\Phi$ over $\Gamma''$ and obtain
\begin{align*}
\int_{\Gamma''}\Phi\dir&=\int_{2A}^{\infty}\int_{B''a}^{\infty}\int_{(1-\varepsilon'')x_1}^{(1+\varepsilon'')x_1}\frac{1}{\log(a^{-1}\,x_1^2)\,(a^{-1}\,x_1^2)^2}\di x_2\di x_1\frac{\di a}{a}\\
&= C\,\int_{2A}^{\infty}\int_{B''a}^{\infty}\frac{x_1}{(a^{-1}\,x_1^2)^2\,\log(a^{-1}\,x_1^2)}\di x_1\frac{\di a}{a}\\
&= C\,\int_{2A}^{\infty}\int_{(B'')^2a}^{\infty}\frac{\di u}{u^2\,\log u}{\di a}\\
&\geq C\,\int_{2A}^{\infty}\frac{1}{a\,\log a}\di a\\
&=\infty\,.
\end{align*} 
Given $(x_1,x_2,a)\in \Gamma''$ and $(y_1,y_2,b),\,(z_1,z_2,c)\in E$, we have
\begin{align*}
&(x_1,x_2,a)\cdot (y_1,y_2,b)^{-1}\cdot(z_1,z_2,c)\\
=&\big(x_1+ab^{-1}(z_1-y_1),x_2+ab^{-1}(z_2-y_2),ab^{-1}c    \big)\,,
\end{align*}
where $ab^{-1}c>2A/2=A$, and
\begin{align*}
x_1+ab^{-1}(z_1-y_1)&>B''a-ab^{-1} > B'' a/2 >B'\,ab^{-1}c\,,
\end{align*}
for $B''$ sufficiently large. Moreover,
\begin{align*}
x_2+ab^{-1}(z_2-y_2)&>x_1(1-\varepsilon'')-ab^{-1}/4\\
&= x_1(1-\varepsilon')+(\varepsilon' -\varepsilon'')x_1 -ab^{-1}/4\\
&> x_1(1-\varepsilon')+(\varepsilon' -\varepsilon'')B''a -ab^{-1}/4\\
&> x_1(1-\varepsilon')+ \big[(\varepsilon' -\varepsilon'')B''-1/4\big]\,ab^{-1}\\
&>\big[{x_1+ab^{-1}(z_1-y_1)}\big]\,(1-\varepsilon')\,,
\end{align*}
if $\varepsilon''<\varepsilon'$ and $B''$ is sufficiently large. In the same way, we can achieve
\begin{align*}
x_2+ab^{-1}(z_2-y_2)&< \big[x_1+ab^{-1}(z_1-y_1)\big]\,(1+\varepsilon')\,.
\end{align*}
Thus the point $(x_1,x_2,a)\cdot (y_1,y_2,b)^{-1}\cdot(z_1,z_2,c)$ is in $\Gamma'$, proving (\ref{inclusion}). 
\end{proof}

\begin{theo}\label{S1}
The operators $S_1$ and $S_2$ are not bounded from $H^1$ to $L^1$.
\end{theo}
\begin{proof}
By symmetry, it is enough to treat the case of $S_1$. We shall construct an atom $a$ such that $S_1a$ does not belong to $L^1$. Let $R$ be the parallelepiped $\big[-\nep^{2}\log 2/2,\nep^{2}\log 2/2 \big]\times \big[-\nep^{2}\log 2/2,\nep^{2}\log 2/2 \big]\times [1/2,2]$; it is easy to check that $R$ is a \CZ set centered at the identity. Now let $E$ be the parallelepiped defined in Lemma \ref{estimateX2k1}, and consider the right translate $E^{\sigma}$ of $E$ by the point ${\rm{exp}}(\sigma\,X_2)=(0,\sigma,1)$ for some $\sigma>0$, i.e.,
\begin{align*}
E^{\sigma}&=E\cdot(0,\sigma,1)=\{(y_1,y_2+b\,\sigma,b):~( y_1,y_2,b)\in E \}\\
&\subset (-1/2,1/2)\times \big(-1/4+\sigma ,2\,\sigma  \big)\times (1,2)\,.
\end{align*}
With $\sigma=1/4$, $E$ and $ E^{\sigma}$ are disjoint and contained in $R$.

Let us consider the function $a=\rho(R)^{-1}\,\big({\bf{1}}_{E}- {\bf{1}}_{E^{\sigma}} \big)$. It is obvious that $a$ is supported in the \CZ set $R$ and  $\|a\|_{\infty}\leq \rho(R)^{-1}$. Moreover $\int a\dir=0$ and so $a$ is an atom. We now compute $S_1a$ outside the support of $a$. For all $x\notin\overline{ E\cup E^{\sigma}}$
\begin{align*}
S_1a(x)&=\int S_1(x,y)\,a(y)\dir(y)\\
&=\rho(R)^{-1}\int_{E} S_1(x,y)\dir(y)-\rho(R)^{-1}\int_{E^{\sigma}} S_1(x,y)\dir(y)\,.
\end{align*} 
Changing variable $y=v\cdot (0,\sigma,1) $ in the last integral, this transforms into
\begin{align*}
&\rho(R)^{-1}\int_{E} S_1(x,y)\dir(y)-\rho(R)^{-1}\int_{E} S_1(x,v\cdot (0,\sigma,1))\dir(v)\\
&=\rho(R)^{-1}\int_{E} \big[S_1(x,y)-S_1\big(x,y\cdot (0,\sigma,1)\big)\big]\dir(y)\,.
\end{align*} 
By (\ref{kernelSi}) we know that 
\begin{align*}
S_1(x,y)-S_1\big(x,y\cdot (0,\sigma,1)\big)&=\delta(x)\,\big(-k_1(x^{-1}y)+k_1(x^{-1}y\, {\rm{exp}}(\sigma\,X_2))  \big)\\
&=\delta(x)\,\sigma\,\frac{d}{dt}\Big\lvert_{ t=\tau(x,y) }k_1\big(x^{-1}y\, {\rm{exp}}(t\,X_2)\big)\\
&= \delta(x)\,\sigma\,X_2k_1\big(x^{-1}y\, {\rm{exp}}(\tau(x,y)\,X_2)\big)\,,
\end{align*}
for some $\tau(x,y)$ in $(0,\sigma)$. It follows that for all $x\notin\overline{ E\cup E^{\sigma}}$
\begin{align}\label{S1a(x)}
S_1a(x)&=\rho(R)^{-1}\,\sigma\,\delta(x)\,\int_{E}X_2k_1\big(x^{-1}y\, {\rm{exp}}(\tau(x,y)\,X_2)\big)\dir(y)\,.
\end{align}
To prove that $S_1a$ is not in $L^1$, we integrate $|S_1a|$ in the region $E\,(\Gamma'')^{-1}$, where $\Gamma''$ is the set which appears in Lemma \ref{estimateX2k1}. It is easy to check that $E\,(\Gamma'')^{-1}$ is disjoint with $\overline{ E\cup E^{\sigma}}$, so that we can apply (\ref{S1a(x)}) in the region $E\,(\Gamma'')^{-1}$ and obtain
\begin{align*}
&\int_{E\,(\Gamma'')^{-1}}\big|S_1a(x)\big|\dir(x)\\
&=\rho(R)^{-1}\,\sigma\,\int_{E\,(\Gamma'')^{-1}}\delta(x)\,\Big|\int_{E}X_2k_1\big(x^{-1}y\, {\rm{exp}}(\tau(x,y)\,X_2)\big)\dir(y)\Big|\dir(x)\\
&=\rho(R)^{-1}\,\sigma\,\int_{\Gamma''\,E^{-1}}\Big|\int_{E}X_2k_1\big(xy\, {\rm{exp}}(\tau(x^{-1},y)\,X_2)\big)\dir(y)\Big|\dir(x)\,.
\end{align*}
If $x\in \Gamma''\,E^{-1}$ and $y\in E$, then $xy\in\Gamma'$, in view of (\ref{inclusion}). Since $0<\tau(x^{-1},y)<\sigma=1/4$, by Lemma \ref{estimateX2k1} the point $xy\, {\rm{exp}}(\tau(x^{-1},y)\,X_2)$ is in $\Gamma$ and
$$X_2k_1\big(xy\, {\rm{exp}}(\tau(x^{-1},y)\,X_2)\big)\geq C\,\Phi(xy\, {\rm{exp}}(\tau(x^{-1},y)\,X_2))= C\,\Phi(xy)\,.$$
Hence, applying Fubini's theorem and using $w=xy$ instead of $x$, we get 
\begin{align*}
\int_{E\,(\Gamma'')^{-1}}\big|S_1a(x)\big|\dir(x)&\geq C\,\rho(R)^{-1}\,\sigma\,\int_{\Gamma''\,E^{-1}}\,\int_{E}\Phi(xy)\dir(y)\dir(x)\\
&=C\,\rho(R)^{-1}\,\sigma\,\int_{E}\dir(y)\int_{\Gamma''\,E^{-1}y}\Phi(w)\dir(w)\\
&\geq C\,\rho(R)^{-1}\,\sigma\,\int_{E}\dir(y)\int_{\Gamma''}\Phi(w)\dir(w)\,.
\end{align*}
Lemma \ref{estimateX2k1} (iii) implies that this integral diverges.
\end{proof}

\section{Unboundedness of $S_0$}\label{Szero}
To prove that the operator $S_0$ is not bounded from $H^1$ to $L^1$, we use the same idea as in the previous section. The only difference is that we consider now the derivative $X_0k_0$ in a slightly different region.

We first compute the derivative of the expression (\ref{k0}) for $k_0$ along the vector field $X_0$:
\begin{align}\label{X0k0}
&X_0k_0(x_1,x_2,a)=\nonumber\\
&=\frac{1}{2\,\pi^{2}}\,\frac{a^{-1}}{r\,\sinh r}+\frac{1}{2\,\pi^{2}}\,\,\frac{1-a^{-2}-a^{-2}(x_1^2+x_2^2)}{2}\,\frac{\sinh r+r\,\cosh r}{r^2\,\sinh^3 r}-\nonumber\\
&-\frac{1}{2\,\pi^{2}}\,\big[a^{-2}+a^{-2}(x_1^2+x_2^2)\big]\,\frac{\sinh r+r\,\cosh r}{r^2\,\sinh^3 r}+\nonumber\\
&+\frac{1}{2\,\pi^{2}}\,\,\frac{-1+a^{-2}+a^{-2}(x_1^2+x_2^2)}{2} \,\frac{a-a^{-1}-a^{-1}(x_1^2+x_2^2)}{2\,\sinh r}\times \nonumber\\
&\times\Big[\,\frac{(2\cosh r+r\,\sinh r )r^2\,\sinh^3 r}{r^4\sinh^6 r}-\nonumber\\
&-\frac{(\sinh r+r\,\cosh r)(2r\,\sinh^3 r+3\,r^2\,\sinh^2 r\,\cosh r)}{r^4\sinh^6 r}\Big]\nonumber\\
&=\frac{1}{2\,\pi^{2}}\,\frac{a^{-1}}{r\,\sinh r}+\frac{1}{2\,\pi^{2}}\,\,\frac{1-3\,a^{-2}-3\,a^{-2}(x_1^2+x_2^2)}{2}\,\frac{\sinh r+r\,\cosh r}{r^2\,\sinh^3 r}+\nonumber\\
&+\frac{1}{2\,\pi^{2}}\,a^{-1}\,\frac{\big[a-a^{-1}-a^{-1}(x_1^2+x_2^2)\big]^2}{4}\times\nonumber\\
&\times\frac{2\,r^2\cosh^2+r^2+2\,\sinh^2 r+3\,r\,\sinh r\,\cosh r}{r^3\sinh^5 r}\,.
\end{align}
\begin{lemm}\label{estimateX0k0}
There exist two regions $\Omega'\subset  \Omega$, a positive continuous function $\Psi $ on $\Omega$ and a positive constant $C$ such that 
\begin{itemize}
\item[(i)] $X_0k_0\geq C\,\Psi $ in ${\Omega}$;
\item[(ii)] for any $(x_1,x_2,a)$ in $\Omega$ and $\tau$ in $[0,1]$, the point $(x_1,x_2,a)\cdot(0,0,\nep^{\tau})$ is in ${\Omega}$ and
$$\Psi\big((x_1,x_2,a)\cdot(0,0,\nep^{\tau})\big)\geq C\,\Psi(x_1,x_2,a)\,;$$
\item[(iii)] $\int_{\Omega'}\Psi \dir=\infty$\,.
\end{itemize}
Let $F$ be the parallelepiped $(-1/{16},1/{16})\times(-1/{16},1/{16}) \times (1,\sqrt{2})$. Then 
\begin{equation}\label{inclusionF}
\Omega'\cdot F^{-1}\cdot F\subseteq \Omega\,.
\end{equation}
\end{lemm}
\begin{proof}
Let $A>1$ be a constant to be chosen later and define
\begin{align}\label{omega}
\Omega&=\{(x_1,x_2,a)\in G:~ x_1^2+x_2^2< a^2/4,\,a>A  \}\,,\nonumber\\
\Omega'&=\{(x_1,x_2,a)\in G:~ x_1^2+x_2^2< a^2/{64},\,a>\sqrt{2}A  \}\,.
\end{align} 
For all $(x_1,x_2,a)$ in ${\Omega}$
$$\frac{a}{2}<\cosh r(x_1,x_2,a) <C\,a\,.$$
For $A$ sufficiently large, $r(x_1,x_2,a)>1$ here, and, since $\nep^{r}\leq 2\cosh r\leq C\,a$, we have $r  \leq C\,\log a$.

It is easy to show that in the region ${\Omega}$ all the summands which appear in the last expression in (\ref{X0k0}) are positive, so that for all $(x_1,x_2,a)$ in ${\Omega}$
$$X_0k_0(x_1,x_2,a)\geq C\,\frac{a^{-1}}{r\,\sinh r}\geq \frac{C}{a^2\,\log a}\,.$$
We define 
$$\Psi(x_1,x_2,a)=\frac{1}{a^2\,\log a}\,.$$ 
The condition (i) is satisfied. 

Let $(x_1,x_2,a)\in \Omega$ and $\tau\in [0,1]$. It is easy to check that the point $(x_1,x_2,a)\cdot(0,0,\nep^{\tau})=(x_1, x_2,a\,\nep^{\tau})$ is in ${\Omega}$. Moreover,
$$\Psi \big((x_1,x_2,a)\cdot(0,0,\nep^{\tau})\big)= \,\frac{1}{a^2\,\nep^{2\tau}\,\log(a\,\nep^{\tau})}   \geq C\,\frac{1}{a^2\,\log a}=C\,\Psi(x_1,x_2,a)\,,$$
as claimed in (ii). To prove (iii), we integrate $\Psi$ over $\Omega'$ and obtain
\begin{align*}
\int_{\Omega'}\Psi\dir&=\int_{\sqrt{2}A}^{\infty}\frac{1}{a^2\log a}\int\int_{x_1^2+x_2^2\leq a^2/{64}}\di x_1\di x_2\frac{\di a}{a}\\
&=C\,\int_{\sqrt{2}A}^{\infty}\frac{1}{a\log a}{\di a}\\
&=\infty\,.
\end{align*} 
Given $(x_1,x_2,a)\in \Omega'$ and $(y_1,y_2,b),\,(z_1,z_2,c)\in F$, we have that 
\begin{align*}
&(x_1,x_2,a)\cdot (y_1,y_2,b)^{-1}\cdot(z_1,z_2,c)\\
=&\big(x_1+ab^{-1}(z_1-y_1),x_2+ab^{-1}(z_2-y_2),ab^{-1}c    \big)\,,
\end{align*} 
where $ab^{-1}c>\sqrt{2}A/{\sqrt 2}=A$, and
\begin{align*}
&\big[x_1+ab^{-1}(z_1-y_1)]^2+\big[x_2+ab^{-1}(z_2-y_2)]^2\\
&< \big(|x_1|+a/8)^2+\big(|x_2|+a/8)^2\\
&< 2\,\big(1/8+1/8)^2\,a^2\\
&<(ab^{-1}c)^2/4\,.
\end{align*}
Thus $(x_1,x_2,a)\cdot (y_1,y_2,b)^{-1}\cdot(z_1,z_2,c)\in \Omega$, and (\ref{inclusionF}) is proved. 
\end{proof}

\begin{theo}\label{S0}
The operator $S_0$ is not bounded from $H^1$ to $L^1$.
\end{theo}
\begin{proof}
Following closely the proof of Theorem \ref{S1}, we shall construct an atom $a$ such that $S_0a$ does not belong to $L^1$. With $R$ as in that proof, we let $F$ be the parallelepiped defined in Lemma \ref{estimateX0k0} and consider the right translate $F^{\sigma}$ of $F$ by the point ${\rm{exp}}(\sigma\,X_0)=(0,0,\nep^{\sigma})$, i.e.,
\begin{align*}
F^{\sigma}&=F\cdot(0,0,\nep^{\sigma})=\{(y_1,y_2,a\nep^{\sigma}):~( y_1,y_2,b)\in F \}\\
&= \big(-1/{16},1/{16})\times \big(-1/{16},1/{16})\times (\nep^{\sigma},\nep^{\sigma}\sqrt{2})\,.
\end{align*}
With $\sigma=(\log 2)/2$, $F$ and $ F^{\sigma}$ are disjoint and contained in $R$.

Let us consider the atom $a=\rho(R)^{-1}\,\big({\bf{1}}_{F}- {\bf{1}}_{F^{\sigma}} \big)$. We compute $S_0a$ outside the support of $a$. For all $x\notin\overline{ F\cup  F^{\sigma}}$
\begin{align*}
S_0a(x)&=\rho(R)^{-1}\int_{F} S_0(x,y)\dir(y)-\rho(R)^{-1}\int_{F^{\sigma}} S_0(x,y)\dir(y)\\
&=\rho(R)^{-1}\int_{F} \big[S_0(x,y)-S_0\big(x,y\cdot (0,0,\nep^{\sigma})\big)\big]\dir(y)\,.
\end{align*} 
By (\ref{kernelSi}) we know that 
\begin{align*}
S_0(x,y)-S_0\big(x,y\cdot (0,0,\nep^{\sigma})\big)&=\delta(x)\,\big(-k_0(x^{-1}y)+k_0(x^{-1}y\, {\rm{exp}}(\sigma\,X_0))  \big)\\
&=\delta(x)\,\sigma\,\frac{d}{dt}\Big\lvert_{ t=\tau(x,y) }k_0\big(x^{-1}y\, {\rm{exp}}(t\,X_0)\big)\\
&= \delta(x)\,\sigma\,X_0k_0\big(x^{-1}y\, {\rm{exp}}(\tau(x,y)\,X_0)\big)\,,
\end{align*}
for some $\tau(x,y)$ in $(0,\sigma)$. It follows that for all $x\notin\overline{ F\cup F^{\sigma}}$
\begin{align}\label{S0a(x)}
S_0a(x)&=\rho(R)^{-1}\,\sigma\,\delta(x)\,\int_{F}X_0k_0\big(x^{-1}y\, {\rm{exp}}(\tau(x,y)\,X_0)\big)\dir(y)\,.
\end{align}
To prove that $S_0a$ is not in $L^1$, we integrate $S_0a$ in the region $F\,(\Omega')^{-1}$. It is easy to verify that $F\,(\Omega')^{-1}$ is disjoint with $\overline{ F\cup F^{\sigma}}$, so that we can apply (\ref{S0a(x)}) and obtain
\begin{align*}
&\int_{F\,(\Omega')^{-1}}\big|S_0a(x)\big|\dir(x)=\\
&=\rho(R)^{-1}\,\sigma\,\int_{F\,(\Omega')^{-1}}\delta(x)\,\Big|\int_{F}X_0k_0\big(x^{-1}y\, {\rm{exp}}(\tau(x,y)\,X_0)\big)\dir(y)\Big|\dir(x)\\
&=\rho(R)^{-1}\,\sigma\,\int_{\Omega'\,F^{-1}}\Big|\int_{F}X_0k_0\big(xy\, {\rm{exp}}(\tau(x^{-1},y)\,X_0)\big)\dir(y)\Big|\dir(x)\,.
\end{align*}
If $x\in \Omega'\,F^{-1}$ and $y\in F$, then $xy\in\Omega$, in view of (\ref{inclusionF}). Since $0<\tau(x^{-1},y)<\sigma<1$, by Lemma \ref{estimateX0k0}(ii) the point $xy\, {\rm{exp}}(\tau(x^{-1},y)\,X_0)$ is in ${\Omega }$ and
$$X_0k_0\big(xy\, {\rm{exp}}(\tau(x^{-1},y)\,X_0)\big)\geq C\,\Psi(xy\, {\rm{exp}}(\tau(x^{-1},y)\,X_0))\geq C\,\Psi(xy)\,.$$
As in the proof of Theorem \ref{S1}, we get
\begin{align*}
\int_{F\,(\Omega')^{-1}}\big|S_0a(x)\big|\dir(x)&\geq C\,\rho(R)^{-1}\,\sigma\,\int_{\Omega'\,F^{-1}}\,\int_{F}\Psi(xy)\dir(y)\dir(x)\\
&=C\,\rho(R)^{-1}\,\sigma\,\int_{F}\dir(y)\int_{\Omega'\,F^{-1}y}\Psi(w)\dir(w)\\
&\geq C\,\rho(R)^{-1}\,\sigma\,\int_{F}\dir(y)\int_{\Omega'}\Psi(w)\dir(w)\,.
\end{align*}
Lemma \ref{estimateX0k0} (iii) implies that the last integral diverges.
\end{proof}


\section{The local parts of $T_{ij}$,\,$S_{ij}$ and $R_{ij}$}\label{local}
In this section, we study the local parts of the kernels of the second-order Riesz transforms. We shall prove that they behave like standard \CZ kernels in $\RR^3$ and deduce that they correspond to operators which are bounded from $H^1$ to $L^1$. 

Let $\Psi$ be a function in $C^{\infty}_c(S)$ such that $0\leq \Psi\leq 1$, $\Psi$ is supported in the ball $B_2$ of radius $2$ and $\Psi=1$ on the ball $B_{1}$. Define
$$g_{ij}^0= g_{ij}\,\Psi\qquad {\rm{and}} \qquad    g_{ij}^{\infty}=g_{ij}\,(1-\Psi)\,,$$
$$k_{ij}^0= k_{ij}\,\Psi\qquad {\rm{and}} \qquad    k_{ij}^{\infty}=k_{ij}\,(1-\Psi)\,,$$
$$\ell_{ij}^0= \ell_{ij}\,\Psi\qquad {\rm{and}} \qquad    \ell_{ij}^{\infty}=\ell_{ij}\,(1-\Psi)\,,$$
and let $T_{ij}^0$, $T_{ij}^{\infty}$, $R_{ij}^0$, $R_{ij}^{\infty}$, $S_{ij}^0$ and $S_{ij}^{\infty}$ be the corresponding convolution operators. We shall prove that the operators $T_{ij}^0$, $R_{ij}^0$ and $S_{ij}^0$ are bounded from $H^1$ to $L^1$. To do so, we use the following lemma.

\begin{lemm}\label{CZest}
Let $T$ be a convolution operator which is bounded on $L^2$. Suppose that its kernel $k$ is a distribution supported in the ball $B_2$ and given by a function in $B_2\setminus \{e \}$. Define 
\begin{equation}\label{beta}
\beta\big( (x_1,x_2,s),\,(y_1,y_2,t) \big)=\delta(y_1,y_2,\nep^{t})\,k\big( (y_1,y_2,\nep^t)^{-1}\cdot_G(x_1,x_2,\nep^s) \big)
\end{equation}
for any $(x_1,x_2,s)\neq (y_1,y_2,t)\in\RR^3$. If $\beta$ satisfies the standard estimate
\begin{equation}\label{CZestimates}
|\beta({\bf{x}},{\bf{y}})|+|{\bf{x}}-{\bf{y}}|\big[ |\nabla_{{\bf{x}}}\beta({\bf{x}},{\bf{y}})|+   |\nabla_{{\bf{y}}}\beta({\bf{x}},{\bf{y}})|\big]\leq C\,|{\bf{x}}-{\bf{y}}|^{-3}\,,
\end{equation}
for $|{\bf{y}}|<2A_0 ,\,{\bf{x}}\neq {\bf{y}}$, where $A_0$ is a suitable constant, then $T$ is bounded from $H^1$ to $L^1$.
\end{lemm}
\begin{proof}
We first verify that the operator $T$ is of weak type $1$. Via a standard Calder\'on-Zygmund decomposition argument, the $L^2$-boundedness of $T$ and the estimate (\ref{CZestimates}) imply that for any $f\in L^1(B_1)$
$$\rho\big(\{ x\in G:|Tf(x)|>t \}  \big)\leq \frac{C}{t}\,\|f\|_1\qquad \forall t>0\,.$$
There exists a sequence of balls $B_j$, centered at points $x_j$ and of radius $1$, such that $G=\bigcup_jB_j$ and each point of $G$ belongs to at most $n$ of the balls $B_j$ (see \cite[Lemma 8]{GQS}). From the left-invariance of the operator $T$, the right-invariance of the measure and a simple application of a partition of unity $(\psi_j)_j$ such that ${\rm{supp}}\,\psi_j \subseteq \overline{B}_j$, we may deduce that for any $f\in L^1$ and $t>0$
\begin{align*}
\rho\big(\{ x\in G:|Tf(x)|>t \}  \big)&\leq \rho\big(\{ x\in G:\sum_j|T(\psi_jf)(x)|>t \}  \big)\\
&\leq    \frac{C}{t}\,\sum_j\|\psi_jf\|_1\leq\frac{C}{t}\,\| f\|_1 \,.
\end{align*}
The inequalities above follow by a standard argument (see \cite[Lemma 7]{GQS} for the details). Thus, $T$ is of weak type $1$. As in the proof of Theorem \ref{Teolim}, the lemma will follow if we show that there exists a constant $C$ such that $\|Tb\|_1\leq C$ for any atom $b$. 

Any atom $b$ can be transformed by an appropriate left-translation into an atom $a$ supported in a \CZ set centered at the identity, and $\|Tb\|_1=\|Ta\|_1$ by the left-invariance of $T$. Thus, it suffices to consider an atom $a$ supported in a \CZ set $R=[-L/2,L/2]\times [-L/2,L/2]\times [\nep^{-r},\nep^{r}]$ centered at the identity. Recall that the dilated set $R^*$ is defined by $\{x\in G:~d(x,R)<r\}$. Since $T$ is bounded on $L^2$, 
\begin{equation}\label{1}
\|Ta\|_{L^1(R^*)}\leq \rho(R^*)^{1/2}\,|\!|\!| T |\!|\!|_{L^2\to L^2}\,\|a\|_2\leq C\,.
\end{equation}
Note that ${\rm{supp}}(Ta)\subseteq R\cdot B_2\subseteq \{x\in S:~d(x,R)<2  \}$. 

If $r>2$, then ${\rm{supp}}(Ta)\subseteq R\cdot B_2\subseteq R^*$, so that $\|Ta\|_{1}=\|   Ta\|_{L^1(R^*)}\leq C$. 

Suppose now that $r<2$. Since $r(x_1,x_2,a)\sim  |(x_1,x_2,\log a)|$ near the identity, there exists an absolute constant $A_0$ such that 
$$R\subset \{(x_1,x_2,a):~~ |(x_1,x_2,\log a)|<A_0r \}=B\,.$$ 
Notice that $B$ corresponds to a euclidean ball in $\RR^3$. Since $r<2$, $\rho(B)\sim  \rho(R)\sim  r^3$ and, by arguing as in (\ref{1}), we obtain that $\|Ta\|_{L^1(2B)}\leq C\,.$

It remains to estimate the $L^1$-norm of $Ta$ outside $2B$. Since the integral kernel of $T$ is $\beta$, we get 
\begin{align*}
\int_{(2B)^c}|Ta(x)|\dir(x)&=\int_{(2B)^c}\Big|\int_{B}a(y)\,\beta(x,y)\dir (y)\Big|\dir(x)\\
&=\int_{(2B)^c}\Big|\int_{B}a(y)\big[ \beta(x,y)-\beta(x,{\bf{0}})  \big]\dir(y)\Big|\dir(x)\,,
\end{align*}
which is bounded by
\begin{align*}
\|a\|_{\infty}\,\int_{(2B)^c}\int_{B}&\big|\beta\big((x_1,x_2,a),(y_1,y_2,b)\big)-\beta\big((x_1,x_2,a),{\bf{0}}\big)  \big|\\
&\dir(y_1,y_2,b) \dir(x_1,x_2,a)\,.
\end{align*}
By changing variables $a=\nep^s$ and $b=\nep^t$, we obtain 
\begin{align*}
&\|Ta\|_{L^1(2B)^c}\leq  \|a\|_{\infty}\times \\
&\times  \int_{ |(x_1,x_2,s) |>2A_0r        }\int_{|(y_1,y_2,t)  |<A_0r}\big|\beta\big((x_1,x_2,s),(y_1,y_2,t)\big)-\beta\big((x_1,x_2,s),{\bf{0}}\big)  \big|\\
&\phantom{ \int_{ |(x_1,x_2,s) |>2A_0r        }\int_{|(y_1,y_2,t)  |<A_0r}~} \di y_1\di y_2\di t \di x_1\di x_2\di s\,.
\end{align*}
If $|{\bf{x}}|>2A_0r$ and $|{\bf{y}}|<A_0r$, by (\ref{CZestimates}) we get
$$\big|\beta\big({\bf{x}},{\bf{y}}\big)-\beta\big({\bf{x}},0\big)  \big|\leq \sup_{|{\bf{y}}'|<A_0r    }|\nabla_{{\bf{y}}}\beta({\bf{x}},{\bf{y}}')|\,|{\bf{y}}|\leq C\,r\,|{\bf{x}}|^{-4}\,,$$
so that
\begin{align*}
\|Ta\|_{L^1(2B)^c}&\leq  \|a\|_{\infty}\,\int_{|{\bf{x}}|>2A_0r}\int_{|{\bf{y}}|<A_0r}\big|\beta\big({\bf{x}},{\bf{y}}\big)-\beta\big({\bf{x}},0\big)  \big|\di {\bf{y}} \di {\bf{x}}\\
&\leq C\,r^{-3}\, \int_{|{\bf{x}}|>2A_0r}\int_{|{\bf{y}}|<A_0r}        |{\bf{x}}|^{-4}\,r\di {\bf{y}}\di {\bf{x}}\\
&\leq C\,.
\end{align*}
\end{proof}
\begin{prop}
The operators $T_{ij}^0$, $R_{ij}^0$ and $S_{ij}^0$ are bounded from $H^1$ to $L^1$.
\end{prop}
\begin{proof}
It is enough to apply Lemma \ref{CZest} to  the operators $T_{ij}^0$, $R_{ij}^0$ and $S_{ij}^0$. By \cite[Theorem 12]{GS1} they are bounded on $L^2$ and their kernels $g_{ij}^0$, $k_{ij}^0$, $\ell_{ij}^0$  are supported in the ball $B_2$. Let $k$ denote one of the kernel $g_{ij}^0$, $k_{ij}^0$, $\ell_{ij}^0$. We must show that the function $\beta$, given by 
$$\beta\big( (x_1,x_2,s),\,(y_1,y_2,t) \big)=\nep^{-2t}\,k\big(\nep^{-t}(x_1-y_1),\nep^{-t}(x_2-y_2),\nep^{s-t}   \big)\,,$$
satisfies (\ref{CZestimates}). By means of some elementary Taylor expansions in the variable ${\bf{x}}=(x_1,x_2,\log a)\in \RR^3$, one finds that near $e$
$$W(x_1,x_2,a)=\frac{1}{4\pi }\,\frac{1}{|\bf{x}|}\,(1+h_1+h_2+\ldots)\,,$$
where each $h_j$ is a function of ${\bf{x}}$ which is homogeneous of degree $j$ and smooth away from $0$, and the series $\sum_jh_j$ converges near $0$. Termwise differentiation is possible, and we let $\partial^{\alpha}$ denotes a differentiation operator with respect to $(x_1,x_2,\log a)$, of order $\alpha$. Then
$$\partial^{\alpha}W=\tilde{h}_{-1-|\alpha|}+\tilde{h}_{|\alpha|}+\ldots\,,$$
with similar smooth homogeneous functions $\tilde{h}_j$. This implies
$$|\partial^{\alpha}W({\bf{x}})|\leq C\,|{\bf{x}}|^{-1-|\alpha|}$$
for small $\bf{x}$, and thus
\begin{equation}\label{nablak}
|k(x_1,x_2,a)|+|(x_1,x_2,\log a)  |\,|\nabla k(x_1,x_2,a)|\leq C\, |(x_1,x_2,\log a)  |^{-3}\end{equation}
in $B_2\setminus \{ e \}$, where $\nabla$ denotes the gradient with respect to $(x_1,x_2,\log a)$.

If $(y_1,y_2,t)$ is near the origin, we have
\begin{align*} 
\big|\beta\big( (x_1,x_2,s),\,(y_1,y_2,t) \big)\big|&\leq C\, \big|k\big(\nep^{-t}(x_1-y_1),\nep^{-t}(x_2-y_2),\nep^{s-t}   \big)\big|\\
&\leq C\, |(\nep^{-t}(x_1-y_1),\nep^{-t}(x_2-y_2),s-t    )  |^{-3}\\
&\leq C\, |(x_1-y_1,x_2-y_2,s-t)  |^{-3}\,,
\end{align*}
and 
\begin{align*}
&\big|\nabla_{(x_1,x_2,s)}\beta\big( (x_1,x_2,s),\,(y_1,y_2,t) \big)\big|+\big|\nabla_{(y_1,y_2,t)}\beta\big( (x_1,x_2,s),\,(y_1,y_2,t) \big)\big|\\
&\leq  C\, \big|k\big(\nep^{-t}(x_1-y_1),\nep^{-t}(x_2-y_2),\nep^{s-t}   \big)\big|+\\
&+ C\,\big|\nabla k\big(\nep^{-t}(x_1-y_1),\nep^{-t}(x_2-y_2),\nep^{s-t}   \big)\big|\\
&\leq C\, |(\nep^{-t}(x_1-y_1),\nep^{-t}(x_2-y_2),s-t    )  |^{-4}\\
&\leq C\, |(x_1-y_1,x_2-y_2,s-t)  |^{-4}\,,
\end{align*}
and the theorem follows.
\end{proof}


\section{Boundedness of $T_{ij}$}\label{Tij}
We shall prove that the operators $T_{ij}=X_i\Delta^{-1}X_j$ are bounded from $H^1$ to $L^1$. Since we already verified the boundedness of their local parts, it remains to consider the global parts. 

In \cite[Lemma 9]{GS1} it is proved that the global parts of the kernels $g_{ij}^{r,\infty}$ of the right-invariant Riesz transforms $T_{ij}^r$ are integrable with respect to the measure $\lambda$. Since (\ref{ZrZell}) implies that $(T_{ij}^r\check{f})^{\lor}=T_{ij}f$ for any $f\in C^{\infty}_c(G)$, we obtain that $g_{ij}^{\infty}=\check{g}_{ij}^{r,\infty}$. Thus $g_{ij}^{\infty}$ is integrable with respect to the measure $\rho$ and the corresponding convolution operator $T_{ij}^{\infty}$ is bounded from $H^1$ to $L^1$. 

\vspace{.5cm}

\section{Unboundedness of $S_{ij}$}\label{Sij}
In this section we prove that the operators $S_{ij}$ are not bounded from $H^1$ to $L^1$. Again it suffices to consider their global parts. To do so, we use the same idea as in Section \ref{Si}, defining an atom whose image under the operator $S_{ij}$ is not integrable far from the support of the atom. 

We will need to estimate some integrals of derivatives of the kernels $k_{ij}$. Notice that it is enough to treat the values of $(i,j)$ listed in the following lemma, since the remaining cases will follow by symmetry.
\begin{lemm}\label{estimateX2kij}
For each pair $(i,j)\in \{(1,1),\,(1,2),\,(1,0),\,(0,1),\,(0,0)  \}$, there exist regions $\Gamma''\subset \Gamma' \subset \Gamma $ in $G$, a positive continuous function $\Phi$ on ${\Gamma}$ and positive constants $C,\tau $ such that 
\begin{itemize}
\item[(i)] $\big|X_2k_{ij}\big|\geq C\,\Phi$ in ${\Gamma}$;
\item[(ii)] for any $(x_1,x_2,a)$ in $\Gamma'$ and $\sigma $ in $[0,\tau]$, the point $(x_1,x_2,a)\cdot(0,\sigma,1)$ is in ${\Gamma}$ and
$$\Phi\big((x_1,x_2,a)\cdot(0,{\sigma},1)\big)=\Phi(x_1,x_2,a)\,;  $$
\item[(iii)] $\int_{\Gamma''}\Phi\dir=\infty\,$.
\end{itemize}
Moreover, there exist constants $0<\delta<1$ and $1<\beta<2$ such that the parallelepiped $E=(0,\delta)\times (-{\delta},0)\times (1,\beta)$ satisfies the condition
\begin{equation}\label{inclusionE}
\Gamma''\cdot E^{-1}\cdot E\subseteq\Gamma'\,.
\end{equation}
\end{lemm}
\begin{proof}
Let us fix a pair $(i,j)$. To simplify the notation we write $k$ for the kernel $k_{ij}$ and drop the indices $i,j$. Because of (\ref{X2kappaij}), there exist constants $\gamma,\,\eta,\,\sigma,\,\theta\in\RR$ and $h,\,\ell,\,m,\,n\in \ZZ^3_+$ such that for $x$ in $\overline{B}_1^c$\begin{align*}
X_2k(x)&=\gamma\,x^h\,\nep^{-2r} +\eta\,x^{\ell}\,\nep^{-3r} +  \sigma\,x^{m}\,\nep^{-3r}+\theta\,x^{n}\,\nep^{-4r} +{Q}(x)=F(x)+Q(x)\,,
\end{align*}
where $F(x)$ is defined by the last equality. Here $\theta\neq 0$, $|h|=0,\,|\ell|=|m|=1,\,|n|=2$. The remainder term $Q(x)$ is as described in Section \ref{kernels}. For large $a$, (\ref{metrica}) implies $\nep^r\sim  a+ a^{-1}(x_1^2+x_2)^2$ and more precisely,
\begin{align*}
\nep^r&=a+a^{-1}(x_1^2+x_2^2)+a^{-1}-\nep^{-r}\\
&=\big( a+a^{-1}(x_1^2+x_2^2) \big)\,\Big(1+\frac{a^{-1}-\nep^{-r}}{a+a^{-1}(x_1^2+x_2^2)}  \Big)\,.
\end{align*}
Inverting the last factor here and expanding, we see that for $x=(x_1,x_2,a)$ with $a$ large
$$\nep^{-pr}= a^p\,(a^2+x_1^2+x_2^2)^{-p}  \big(1+ O(\nep^{-r})  \big)=     a^{p}\,|x|^{-2p}\big(1+ O(\nep^{-r})  \big)\,,
$$
where $|\cdot|$ denotes the euclidean norm in $\RR^3$. Thus, for such $x$ 
\begin{align*}
F(x)&=\frac{ \theta\,x^{n}\,a^{4} +  |x|^2\big[    \gamma\,x^h\,a^{2}\,|x|^2+\eta\,x^{\ell}\,a^{3} +  \sigma\,x^{m}\,a^{3}\big]            }{|x|^8}+E(x)\\
& =\frac{P(x)}{|x|^8}+E(x)\,,
\end{align*}
where $P$ is a polynomial in the variables $x_1,x_2,a$, homogeneous of degree $6$. Further, $E(x)$ is a sum like $F(x)$, but with $\nep^{-(p+1)r}$ instead of $\nep^{-pr}$ in each term. We write $P(x)=\theta\,x^{n}\,a^{4}+|x|^2\,\tilde{P}(x)$, where $\tilde{P}$ is homogeneous of degree $4$. Notice that $P$ is not identically $0$, since the monomial $x^n\,a^4$ cannot equal a product $\theta^{-1}\,|x|^2\,\tilde{P}(x)$. We can thus find $q_1,{q}_2>0$ with $P(q_1,q_2,1)\neq 0$. By continuity and homogeneity, $P(x)\neq 0$ also for $x$ in a narrow cone near the ray in the direction $(q_1,q_2,1)$, in particular for $x$ in the truncated cone
$$\Gamma=\{ (x_1,x_2,a)\in G:~a>A,\,\big|x_1/a-q_1\big|<\varepsilon,\,\big|x_2/a-{q}_2\big|<\varepsilon \}\,,$$
for some small $\varepsilon>0$. With $A$ sufficiently large, this implies that in the region $\Gamma$ the quantities $|E|,\,|Q|$ are much smaller than $|F|$ and so for any $(x_1,x_2,a)\in\Gamma$
$$|X_2k(x_1,x_2,a)|\geq C\,|F(x_1,x_2,a)|\geq C\,\frac{|P(x_1,x_2,a)|}{|(x_1,x_2,a)|^8}\geq C\,a^{-2} \,.$$
Defining $\Phi(x_1,x_2,a)=a^{-2}$ in $\Gamma$, we have proved (i). 

We define
$$\Gamma'=\{ (x_1,x_2,a)\in \Gamma:~\big|x_1/a-q_1\big|<\varepsilon/2,\,\big|x_2/a-{q}_2\big|<\varepsilon/2 \}\,,$$
and
$$\Gamma''=\{ (x_1,x_2,a)\in \Gamma:~a>2A,\,\big|x_1/a-q_1\big|<\varepsilon/4,\,\big|x_2/a-{q}_2\big|<\varepsilon/4 \}\,.$$
Now choose $\tau<\varepsilon/2$ and $\sigma \in [0,\tau]$, and let $(x_1,x_2,a)$ be in $\Gamma'$. Then $(x_1,x_2,a)\cdot (0,\sigma,1)=(x_1,x_2+a\sigma,a)$. We have that $\big|x_1/a-q_1  \big|<\varepsilon/2<\varepsilon$ and
$$\big|(x_2+a\sigma)/a-{q}_2  \big|<\big|x_2/a-{q}_2  \big|+\sigma<\varepsilon/2+\tau <\varepsilon\,.$$
Thus $(x_1,x_2+a\sigma,a)\in {\Gamma}$ and $\Phi\big((x_1,x_2,a)\cdot(0,\sigma,1)\big)=a^{-2}= \Phi(x_1,x_2,a)$.

To prove (iii), it suffices to note that
$$\int_{\Gamma''}a^{-2}     \di x_2\di x_1\frac{\di a}{a}\geq C\,\int_{2A}^{\infty} a^{-2}\,a^2    \frac{\di a}{ a}=\infty\,.$$
Aiming at (\ref{inclusionE}) we take points $(x_1,x_2,a)\in\Gamma''$ and $(y_1,y_2,b),\,(z_1,z_2,c)\in E$ and consider $(x_1,x_2,a)\cdot (y_1,y_2,b)^{-1}\cdot (z_1,z_2,c)=(x_1+ab^{-1}(z_1-y_1),x_2+ab^{-1}(z_2-y_2),ab^{-1}c)$. Obviously, $ab^{-1}c>2A/{\beta}>2A/2=A$. Moreover,
\begin{align*}
\Big|\frac{x_1+ab^{-1}(z_1-y_1)}{ab^{-1}c}-q_1  \Big|&\leq \Big|\frac{x_1}{ab^{-1}c}-q_1  \Big|+\frac{|z_1-y_1|}{c}\\
&\leq \frac{x_1}{a}\,\Big|b/c  -1\Big|+\Big|\frac{x_1}{a}-q_1  \Big|  +\delta\\
&\leq 2q_1\,|\beta-1|+\varepsilon/4+\delta\\
&<\varepsilon/2\,,
\end{align*}
for $\delta$ sufficiently small and $\beta$ sufficiently close to $1$.
In a similar way, we can achieve 
$$\Big| \frac{x_2+ab^{-1}(z_2-y_2)}{ab^{-1}c}-q_2 \Big|<\varepsilon/2\,,$$
so that $(x_1,x_2,a)\cdot (y_1,y_2,b)^{-1}\cdot (z_1,z_2,c)\in \Gamma'$, proving (\ref{inclusionE}). 
\end{proof}

\begin{theo}\label{S10}
The operators $S_{ij}$, for $i,j=0,1,2,$ are not bounded from $H^1$ to $L^1$.
\end{theo}
\begin{proof}
As remarked above, we need only consider the operators $S_{11}$, $\,S_{21}$,$\,S_{10}$,$\,S_{01}$,$\,S_{00}$. 

We argue as in the proof of Theorem \ref{S1}. This time by (\ref{kernelSij}) for all $\sigma>0$ and $x,y\in G$, with $x\neq y$ and $x\neq y\cdot (0,\sigma,1)$ 
\begin{align}\label{diffSij}
S_{ij}(x,y)-S_{ij}\big(x,y\cdot (0,\sigma,1)\big)&= \delta(x)\,\sigma\,X_2k_{ji}\big(x^{-1}y\, {\rm{exp}}(\tau(x,y)\,X_2)\big)\,,
\end{align}
where $0<\tau(x,y)<\sigma$ and $S_{ij}$ denotes the integral kernel of the operator $S_{ij}$.

As in the proof of Theorem \ref{S1} one constructs an atom $a$ such that $\int S_{ij}(\cdot ,y)a(y)\dir (y)$ does not belong to $L^1$: it suffices to apply (\ref{diffSij}) and Lemma \ref{estimateX2kij}. We omit the details.
\end{proof}


\section{Unboundedness of $R_{ij}$}\label{Rij}
In this section we prove that the operators $R_{ij}$ are not bounded from $H^1$ to $L^1$, and it suffices to consider their global parts. The proof of the unboundedness of $R_{ij}^{\infty}$ is different from the previous cases. We shall now construct a sequence of functions in $H^1$ such that their images under the operator $R_{ij}^{\infty}$ lie in $L^1$ but have large $L^1$-norms. To do so, we first analyze the kernels $k_{ij}^{\infty}$.
\begin{lemm}\label{k1k2k3}
For any $i,j=0,1,2,$ there exists a splitting $k_{ij}^{\infty}=k_{ij}^1+k_{ij}^2+k_{ij}^3$ such that
\begin{itemize}
\item[(i)] $k_{ij}^1=k_{ij}^{\infty}\,\chi_{\{(x_1,x_2,a)\in G:~a\leq 1 \}}$ is integrable;
\item[(ii)] $k_{ij}^2$ is supported in the region $\{(x_1,x_2,a)\in B_1^c:~a\geq  1 \}$ and is integrable;
\item[(iii)] $k_{ij}^3$ is supported in the region $\{(x_1,x_2,a)\in B_1^c:~a\geq 1 \}$ and for any $f\in L^1$ 
$$f\ast k_{ij}^3(x_1,x_2,a)=[\psi_a\ast_{\RR^2}h] (x_1,x_2)\qquad \forall (x_1,x_2,a)\in G \,,$$
where $h(x_1,x_2)=\int_{0}^{\infty}f(x_1,x_2,a){\di a}/{a}$, $\psi$ is a continuous function on $\RR^2$ such that $|\psi(x_1,x_2)|\leq C\,(1+|(x_1,x_2)|)^{-3}$ for some $C$, and $\psi_a (x_1,x_2)=a^{-2}\psi(a^{-1}x_1,a^{-1}x_2)$ for $a>0$. 
\end{itemize}
\end{lemm} 
\begin{proof}
We fix a pair $(i,j)$ and drop the indices $i,j$ on the kernels. By (\ref{kappaij}) there exist constants $\alpha,\,\beta$ and $m,\,n\in \ZZ^3$ such that in $B_1^c$
$$
k^{\infty}(x)=\alpha\,x^m\,\nep^{-2r} +\beta\,x^n\,\nep^{-3r} +Q(x)\,,
$$
where $Q$ is integrable, $\beta>0$, $|m|=0$, $|n|=1$ and
\begin{align}\label{mn}
&m_0+2>0\qquad n_0+3>0\nonumber\\      
&m_1+m_2-4< -2 \qquad {\rm{and}} \qquad n_1+n_2-6< -2\,. \end{align}
We define $k^1=k^{\infty}\,\chi_{\{(x_1,x_2,a)\in G:~a\leq 1 \}}$. By (\ref{mn}) and Lemma \ref{integrability and derivative}(ii), $k^1$ is integrable. 

We now consider the region $\{(x_1,x_2,a)\in B_1^c:~a\geq  1 \}$. There we may approximate $\nep^{-r}$ by $1/(2\cosh r)$ and $\cosh r$ by $a\,\big(1+\big|a^{-1}(x_1,x_2)\big|^2\big)$. Estimating the errors, we can write the principal terms in the expression for $k^{\infty}$ above as
\begin{align*}
\alpha\,x^m\,\nep^{-2r} +\beta\,x^n\,\nep^{-3r} &=\alpha\,a^{|m|}\big(a^{-1}x_1\big)^{m_1}\big(a^{-1}x_2\big)^{m_2}\frac{1}{4a^2(1+|a^{-1}(x_1,x_2)|^2)^2}\\
&+ \beta\,a^{|n|}\big(a^{-1}x_1\big)^{n_1}\big(a^{-1}x_2\big)^{n_2}\frac{1}{8a^3(1+|a^{-1}(x_1,x_2)|^2)^3}\\
&+q(x)\\
&=a^{-2}\psi\big(a^{-1}(x_1,x_2)\big)+q(x)\,,
\end{align*}
where 
$$\psi(x_1,x_2)=\alpha \,\frac{x_1^{m_1}x_2^{m_2}}{4(1+|(x_1,x_2)|^2)^2}+\beta\,\frac{x_1^{n_1}x_2^{n_2}}{8(1+|(x_1,x_2)|^2)^3}\,$$
and
$$q(x)=O\big(x^{m}\nep^{-3r}+x^{n}\nep^{-4r}\big)\,.$$ 
By Lemma \ref{integrability and derivative}(i) and (\ref{mn}), $q$ is integrable in the region where $a\geq 1$, and $|\psi(x_1,x_2)|\leq C\,(1+|(x_1,x_2)|)^{-3}$.

Define $k^3(x_1,x_2,a)=a^{-2}\psi\big(a^{-1}(x_1,x_2)\big)\chi_{\{(x_1,x_2,a)\in B_1^c:~a\geq  1 \}}$ and $k^2=k^{\infty}-k^3$. Then 
$$k^2(x_1,x_2,a)=Q(x)+q(x) \,,$$
and so $k^2$  is integrable, which proves (ii).

Given a function $f$ in $L^1$, we obtain
\begin{align*}
&f\ast k^3(x_1,x_2,a)=\\
&=\int_0^{\infty }\int\int_{\RR^2}f(x_1-ab^{-1}y_1,x_2-ab^{-1}y_2,ab^{-1})\,b^{-2}\,\psi(b^{-1}y_1,b^{-1}y_2)\di y_1\di y_2 \di b/b\\
&=\int_0^{\infty }\int\int_{\RR^2}f(x_1-az_1,x_2-az_2,ab^{-1})\,\psi(z_1,z_2)\di z_1\di z_2 \di b/b\\
&=\int\int_{\RR^2}a^{-2}\,\psi(a^{-1}v_1,a^{-1}v_2)   \int_0^{\infty}  f(x_1-v_1,x_2-v_2,c)\di c/c \di v_1\di v_2  \\
&=\big[\psi_a\ast_{\RR^2}h\big](x_1,x_2)\,,
\end{align*}
which proves (iii).
\end{proof}
We remark that in \cite[Section 7]{GS1}, the analog of Lemma \ref{k1k2k3} was proved for the operators $R_{ij}^r=X_i^rX_j^r(\Delta^r)^{-1}$. We could also deduce Lemma \ref{k1k2k3} from that result.

We shall need the following technical lemma, which shows how to construct functions in $H^1(G)$ from functions in $H^1(\RR^2)$.
\begin{lemm}\label{H1}
For any function $h$ in $H^1(\RR^2)$, there exists a function $f$ in $H^1(G)$ such that $\|f\|_{H^1(G)} \leq \|h\|_{H^1(\RR^2)}$ and 
$$h(x_1,x_2)=\int_0^{\infty}f(x_1,x_2,a)\frac{\di a}{a}\,.$$
\end{lemm}
\begin{proof}
Let $h$ be in $H^1(\RR^2)$. Take a decomposition of $h$ as $\sum_j\lambda_jb_j$, where $\lambda_j\in \CC$, $\sum_j|\lambda_j|<\infty$ and $b_j$ are atoms in $\RR^2$. The atom $b_j$ is supported in a square $Q_j$ of side $L_j$, and $\int b_j=0$ and $\|b_j\|_{\infty}\leq L_j^{-2}$. We choose $r_j>0$ such that either $r_j<1$ and $\nep^{2}\,r_j\leq L_j< \nep^{8 }\,r_j$ or $r_j\geq 1$ and $\nep^{2r_j}\leq L_j< \nep^{8r_j}$. Define 
$$a_j(x_1,x_2,a)=\frac{1}{2}\,r_j^{-1}\chi_{[\nep^{-r_j},\nep^{r_j}]}(a)\,b_j(x_1,x_2)\,.$$
The functions $a_j$ are atoms in $G$ supported in the \CZ sets $R_j=Q_j\times[\nep^{-r_j},\nep^{r_j}]$. Now define $f=\sum_j\lambda_ja_j$. It is easy to check that $f$ is $H^1(G)$ and has the required properties.
\end{proof}
We now concentrate on the part of the kernel which is not integrable, i.e., $k_{ij}^3$.
\begin{lemm}\label{k3}
The operator $f\mapsto f\ast k^3_{ij}$ is not bounded from $H^1$ to $L^1$.
\end{lemm}
\begin{proof}
The proof will follow those of \cite[Lemmata 13,\,14]{GQS}. We will define a sequence of functions $h_N$ in the Hardy space $H^1(\RR^2)$ such that ${ \| \psi_a\ast h_N\|_1}/{\|f_N\|_{H^1(\RR^2)}} $ is large. From Lemma \ref{H1}, we then obtain a sequence of functions $f_N$ in $H^1(G)$ such that ${\|f_N\ast k^3_{ij}\|_1}/{\|f_N\|_{H^1}}$ is not uniformly bounded.

Let $\phi$ be a $C^{\infty}$-function in $\RR^2$ supported in $[-1,1]\times [-1,1]$ such that $\int \phi=0$, and $\psi\ast_{\RR^2}\phi(0,0)\neq 0$. Let $L>1$ and let $N$ be the greatest natural number with $N<\log L$. Let $p,\,q$ be large natural numbers to be chosen later. Define
\begin{equation}\label{hN}
h_N=\sum_{n=0}^N\sum_{k\in \ZZ^2,\,|k_i|<(2^{qn}L-1)/p}\pm \phi_{nk}\,,
\end{equation}  
where the signs will be chosen later and $\phi_{nk}(x_1,x_2)=\phi(2^{qn}x_1-pk_1,2^{qn}x_2-pk_2)$. Let $n$ and $k$ be as in the double sum. Since ${\rm{supp}}\,\phi\subset [-1,1]\times [-1,1]$, we conclude
$${\rm{supp}}\,\phi_{nk}\subset [2^{-nq}(pk_1-1),2^{-nq}(pk_1+1) ]\times [2^{-nq}(pk_2-1),2^{-nq}(pk_2+1)]\,.$$ 
It follows that $h_N$ is supported in $[-L,L]\times [-L,L]$. 

{\bf{Claim 1.}} One can choose $p,q$ and $t>0$ independently of $N$ so that for all sign choices in (\ref{hN})
$$\rho(\{(x_1,x_2,a):~|\psi_a\ast h_N(x_1,x_2)|>t \})\geq C\,N\,L^2\,.$$

{\bf{Claim 2.}} The signs in (\ref{hN}) can be chosen so that $\|h_N\|_2\leq C\,\sqrt{N}\,L$.

\smallskip 
{\emph{Proof of Claim 1.}} Since $\psi\ast_{\RR^2}\phi(0,0)\neq 0$, there exists a positive $\delta$ such that $\big|\psi_a\ast_{\RR^2}\phi(x_1,x_2)\big|>\delta$ for $(x_1,x_2,a)$ in a neighbourhood $U$ of $(0,0,1)$ in $\RR^2\times \RR^+=G$. We can take $U$ contained in $[-1,1]\times [-1,1]\times [1/2,2]$. It follows that 
\begin{equation}\label{Unk}
\big|\psi_a\ast_{\RR^2}\phi_{nk}(x_1,x_2)\big|>\delta \qquad{\rm{if~}} (x_1,x_2,a)\in U_{nk}\,,
\end{equation}
where $U_{nk}=\{(x_1,x_2,a):~(2^{qn}x_1-pk_1,2^{qn}x_2-pk_2,2^{qn}a )\in U \}$. The sets $U_{nk}$ are mutually disjoint and $\rho(U_{nk})=2^{-2qn}\rho(U)$. 

Now fix $0\leq m\leq N$, $\ell\in \ZZ^2$ such that $|\ell_i|<(2^{qm}L-1)/p$ and take $(x_1,x_2,a)\in U_{m\ell}$. By (\ref{Unk}), in the sum 
$$\sum_{n=0}^N\sum_{k\in \ZZ^2,\,|k_i|<(2^{qn}L-1)/p}\pm \psi_a\ast \phi_{nk}(x_1,x_2)\,,$$
the term with $n=m,\, k=\ell$ is greater than $\delta$ in absolute value. The other terms are much smaller; more precisely, we can choose $p,q$ such that
\begin{equation}\label{smaller}
\mathop{\sum \sum}_{(n,k)\neq(m,l)}  | \psi_a\ast \phi_{nk}(x_1,x_2)\big|\leq \delta/2\,.
\end{equation}
The proof of (\ref{smaller}) is the same as \cite[Proof of Claim 1, page 277]{GQS}, and we omit it. 

This means that
$$\{(x_1,x_2,a):~|\psi_a\ast h_N(x_1,x_2,a)|>\delta/2  \}\supseteq \bigcup_{n=0}^N\bigcup_{k\in \ZZ^2, |k_i|\leq (2^{qn}L-1)/p}U_{nk}\,.$$
Thus, choosing $t=\delta/2$,
\begin{align*}
\rho(\{(x_1,x_2,a):~|\psi_a\ast h_N(x_1,x_2,a)|>t  \})&\geq \sum_{n=0}^N\sum_{k\in \ZZ^2,\,|k_i|<(2^{qn}L-1)/p}\rho(U_{nk})\\
&\geq C\,\rho(U)\,\sum_{n=0}^N2^{-2qn}(2^{qn}L-1)^2\\
&\geq C\,N\,L^2\,.
\end{align*}

{\emph{Proof of Claim 2.}} This proof follows the idea of \cite[Proof of Claim 2, page 279]{GQS}. On the set of all sign choices in (\ref{hN}), consider the probability measure which makes the signs into independent Bernoulli variables. Denote by ${\bf{E}}$ the corresponding expectation. Then
\begin{align*}
{\bf{E}}|h_N|^2(x)&= \sum_{n,k}|\phi_{nk}(x)|^2\\
&= \sum_{n=0}^N\sum_{k\in \ZZ^2,\,|k_i|<(2^{qn}L-1)/p}|\phi(2^{qn}x_1-pk_1,2^{qn}x_2-pk_2)|^2\,.
\end{align*}
Since $\phi_{nk}$ is supported in $[(pk_1-1)2^{-nq},(pk_1+1)2^{-nq}]\times [(pk_2-1)2^{-nq},(pk_2+1)2^{-nq} ]$ and uniformly bounded, we get
\begin{align*}
&\int {\bf{E}}|h_N(x_1,x_2)|^2\di x_1\di x_2\\
&\leq C\,\sum_{n=0}^N\sum_{k\in \ZZ^2,\,|k_i|<(2^{qn}L-1)/p}\big|[(pk_1-1)2^{-nq},(pk_1+1)2^{-nq}]\times \\
&\phantom{C\,\sum_{n=0}^N\sum_{k\in \ZZ^2,\,|k_i|<(2^{qn}L-1)/p}}\times
 [(pk_2-1)2^{-nq},(pk_2+1)2^{-nq} ] \big|\\
&\leq C\,\sum_{n=0}^N\sum_{k\in \ZZ^2,\,|k_i|<(2^{qn}L-1)/p}2^{-2nq}\\
&\leq C\,N\,L^2\,.
\end{align*}
Thus ${\bf{E}}\|h_N\|_2\leq C\,\sqrt{N}\,L$, and Claim 2 follows. 

\smallskip
If we choose $p,q$ and the signs in (\ref{hN}) as in Claim 1 and 2, the function $h_N$ will be a multiple of a $(1,2)$-atom in $\RR^2$ (see \cite{CW}). Indeed, it is supported in $[-L,L]\times [-L,L]$, with integral zero, and $\|h_N\|_2\leq C\,\sqrt{N}\,L$. In particular, $h_N$ is in $H^1(\RR^2)$ and $\|h_N\|_{H^1(\RR^2)}\leq C\,L^2\,\sqrt{N}$. 

By Lemma \ref{H1}, there exists $f_N$ in $H^1(G)$ such that $\|f_N\|_{H^1(G)} \leq C\,L^2\,\sqrt{N}$ and 
$$h_N(x_1,x_2)=\int_0^{\infty}f_N(x_1,x_2,a)\frac{\di a}{a}\,.$$
Thus by Lemma \ref{k1k2k3} and Claim 1,
\begin{align*}
\|f_N\ast k^3\|_1&=\int_0^{\infty}\int_{\RR^2}|\psi_a\ast h_N(x_1,x_2,a)|\dir(x_1,x_2,a)\\
&\geq t\,\rho(\{(x_1,x_2,a):~|\psi_a\ast h_N(x_1,x_2,a)|>t  \})\\
&\geq C\,N\,L^2\,.
\end{align*}
This shows that ${\|f_N\ast k^3_{ij}\|_1}/{\|f_N\|_{H^1(G)}}$ is not uniformly bounded, proving the lemma.
\end{proof}
\begin{theo}
The operators $R_{ij}^{\infty}$, for $i=0,1,2,$ are not bounded from $H^1$ to $L^1$.
\end{theo}
\begin{proof}
This is a direct consequence of Lemmata \ref{k1k2k3} and \ref{k3}.
\end{proof}

\vspace{1cm}
\noindent
PETER SJ\"OGREN\\ 
Department of Mathematical Sciences\\
G\"oteborg University and\\
Chalmers University of Technology\\
S-412 96 G\"oteborg\\Sweden
\\
peters@math.chalmers.se 
\\
\\
\noindent
MARIA VALLARINO\\
Dipartimento di Matematica e Applicazioni
\\ Universit\`a di Milano-Bicocca\\
Via R.~Cozzi 53\\ 20125 Milano\\ Italy
\\
maria.vallarino@unimib.it

\end{document}